\documentclass[a4paper,12pt]{article}
\usepackage{latexsym, amsmath, amssymb, amsthm}
\newcommand{\Ds}{\displaystyle}
\begin{document}

\newcommand{\lt}{\left }
\newcommand{\rt}{\right }
\renewcommand{\le}{\leqslant}
\renewcommand{\ge}{\geqslant}

\newcommand{\R}{\bf R }
\newcommand{\N}{\bf N}
\newcommand{\Z}{\bf Z}
\newcommand{\C}{\bf Z}

\newcommand{\D}{\Delta}
\renewcommand{\d}{\delta}
\newcommand{\g}{\gamma}
\newcommand{\G}{\Gamma}
\newcommand{\Th}{\varTheta}
\renewcommand{\l}{\lambda}
\newcommand{\eps}{\varepsilon}
\newcommand{\ff}{\varphi}
\renewcommand{\a}{\alpha}
\renewcommand{\b}{\beta}

\newtheorem{thm}{Theorem}
\newtheorem{lm}[thm]{Lemma}
\newtheorem{cor}[thm]{Corollary}
\newtheorem{rem}[thm]{Remark}
\newenvironment{pf}{\par\medskip\noindent{\bf Proof.}}%
{\nopagebreak\hfill$\Box$\par\medskip}

\title{NOISY FIGHTER-BOMBER DUEL}
\author{L.~N.~Positselskaya}
\date{}
\maketitle

\begin{abstract}

\vskip 1mm
  We discuss a duel-type game in which Player~I uses his
resource continuously and Player~II distributes it by discrete
portions.
 Each player knows how much resources he and his opponent have
at every moment of time.
 The solution of the game is given in an explicit form.
\vskip 1,5mm

\end{abstract}

\vskip 1,5mm
\section {Introduction}
 The classical duel is a zero-sum game of two players of
the following type.
 The players have certain resources and use them during a given time
interval with the goal of achieving success.
 Use of the resource $\gamma$ at the moment~$t$ leads to success with
the probability depending on the amount of resource $\gamma$ and
the time $t$ only (it is usually assumed that the probability of
success increases with time).
 As soon as one player reaches the goal he receives his profit,
which is equal to his opponent's loss, and the game ends.
 Various assumptions about the ways the players use their resources
and about the players receiving information about the opponent's
behavior during the game define various kinds of
duels~\cite{kar,kimoverview}.
 Models were considered where the players' resources were discrete
({\it discrete firing\/} duels), infinitely divisible
({\it continuous firing\/} duels), continuous for one of the players
and discrete for the other one ({\it mixed\/} duels, or
{\it fighter-bomber duels\/})~\cite{pdn,pnr}.
 Researchers studied {\it noisy\/} duels~\cite{pdn,fk}, where every
player at a given moment of time had complete information
about his opponent's behavior up to that moment, and {\it silent\/}
duels, where no such information was available.
 At the present time, duels are considered as classical models of
competition~\cite{kar,pmixsin}.

 The author is grateful to Leonid Positselski for his help in
editing this paper and translating it into English.

\section{Posing the problem}

 We consider a zero-sum two-person game of the following structure.
 The players have resources $a\ge 0$ and $m\ge0$ ($m\in \Z$),
respectively.
 The players use their resources during the time interval $[0,1]$ with
the goal of achieving success.
 Player I has an infinitely divisible resource $a\in\R$; $a>0$,
which he uses continuously.
 Player II has a discrete resource $m \in \N$ and uses it by units.
 The effectiveness of the $j$-th player using his resource is described
by the function $P_j(t)$  ($j=1,2$), which defines the probability of
achieving success when using the unit of resource at the moment~$t$.
 The functions $P_j(t)$ are called the {\it effectiveness functions},
they are continuously differentiable and increasing, $P_j(0)=0$,
$P_j(1)=1$, $P_j(t)<1$ for $t<1$.
 Put $p(t)=1-P_1(t)$, $q(t)=1-P_2(t)$, $P(t)=(P_1(t), P_2(t))$.
 The vector-function $P(t)$
is called the {\it effectiveness vector-function}.
 The probability  $G(t,\Delta\gamma)$ of achieving success when using
the amount of resource $\Delta\gamma\ge0$ at the moment~$t$ with
the effectiveness function $P_j(t)$ is described by
the formula~\cite{psi,lk}:
\begin{equation}
 G(t,\Delta\gamma)=1-(1-P_j(t))^{\Delta\gamma},\quad
 \Delta\gamma>0;\quad  G(t,0)=0.
\label{I2}
\end{equation}
 If one of the players achieves success, the game stops.
 The profit of the $j$-th player in the case of his success is equal
to $A_j$.
 Denote by $A$ the pair $(A_1,A_2)$.
 The players' profits are equal to $0$ if none of them has
achieved success or if success has been achieved by both of
them simultaneously.

 Let $\alpha (t)$, $n(t)$ be the players' remaining resources at
the moment $t$.
 Let us call the functions $\alpha (t)$, $n(t)$
the {\it consumption functions\/} of the players.
 They are nonincreasing, and $n(t)$ is a step-function continuous
from the left.
 The pair $(\alpha (t), n(t))$ is called a {\it play\/} of the game.
 Assume that the function $\a(t)$ is continuous in $[0,1]$ and
piecewise continuously differentiable in $(0,1)$.
 Put $\xi(t)=-\a'(t)$ and name $\xi(t)$ the {\it intensity function}.

 A time moment when a player uses his resource is called an {\it action
moment} of this player.
 It is a decreasing point of Player I's consumption function
($\xi(t+0)>0$) or it is a jump point of Player II's consumption
function.
 Let $\eta_i$, $i=1,2,\dots,m$
($0\le\eta_m\le\eta_{m-1}\le\dots\le\eta_1\le 1$) be Player II's
action moments.
 The vector $\eta=(\eta_1,\,\eta_2,\dots\eta_m)$ is called
the {\it vector of action moments}.

 Let us compute the probability $\varphi(\alpha, t_1, t_2)$
of Player I achieving success when using his infinitely divisible
resource according to the consumption function $\alpha(t)$
at the time interval $[t_1,t_2]$, assuming that Player II does
not act during that period of time.
 By \eqref{I2} we have:
\begin{align}
\nonumber
\varphi(\xi, t_1, t_2)&=
\lim_{N\to\infty}\lt(1-\prod_{i=1}^N\lt(1-P_1(\tau_i)\rt)^
{\xi(\tau_i)\D\tau_i}\rt)=\\
\label{phi}
&=1-\exp\int_{t_1}^{t_2}\xi(\tau)\log (1-P_1(\tau))\,d\tau,
\end{align}
where $\xi(t)=-\alpha'(t)$.
 The probability $\varphi(\alpha, t_1, t_2)$ of Player I achieving
success in the segment $[t_1,t_2]$ can be expressed in terms of
the consumption function $\alpha(t)$ as follows:
\begin{align}
\label{phi1}
\varphi(\alpha, t_1, t_2)=
1-\exp\left(-\int_{t_1}^{t_2}\log (1-P_1(\tau))\,d\alpha(\tau)\right).
\end{align}

 Let $K(\alpha,\eta)$ be the mathematical expectation of the profit
received by Player I in the case when Player I uses his resource
according to the consumption function $\alpha(t)$ and Player II
uses his resource at the moments of time $\eta_k$ ($1\le k\le m$).
 It is computed in the following way.
 For $m=0$ we have $K=0$ if $a=0$  and $K=A_1$ if $a>0$.
 For $m\ge 1$ we obtain $K(\alpha,\eta)$ from the recursive formula
\begin{align}
\label{defK}
K(\alpha; \eta_1, \eta_2,\dots, \eta_m)&=
A_1\varphi(\alpha, 0, \eta_m)-A_2(1-\varphi(\alpha, 0, \eta_m))
P_2(\eta_m)+\\
&+(1-P_2(\eta_m))(1-\varphi(\alpha, 0, \eta_m))
K(\alpha_m; \eta_1, \eta_2,\dots, \eta_{m-1})\nonumber,
\end{align}
where
$$
\alpha_m(t)=
  \begin{cases}
    \alpha(\eta_m), &t\in[0,\eta_m);\\
    \alpha(t), &t\in[\eta_m,1].
  \end{cases}
$$

 The game under consideration is called the noisy fighter-bomber duel.
 It is a model of competition in the conditions of complete information.
 In this game every player at a given moment of time has information
about both player's resources up to that moment and continuously
corrects his behavior on the basis of the received information about
the present amount of his opponent's remaining resource.
 Player I's strategy is a function $\xi=u(t,\a,n)$ which determines
the intensity of resource consuming $\xi$ at a moment $t$
in dependence of the current values of players' resources $\a$ and $n$.
 Player II's strategy is a function $\eta_n=v(\a,n)$ assigning
the moment of next action to a pair of players' current resources 
$\a$ and $n$.
 We will define the players' strategies in the segment where
$\alpha(t)n(t)>0$ only and assume that if one of the players has not
exhausted his resource then he consumes it so that the probability
of his success is equal to $1$.
 By the condition $P_j(1)=1$  and formula \eqref{phi1},
this is always possible.
 The payoff function of the game is the function $K(\xi,\eta)$ defined
by the formula \eqref{defK}, where $\xi$ are $\eta$ are the intensity
function of Player I and the action moment vector of Player II
realized during the game.
 Let us denote the game so described by $G_{am}(P, A)$.

\section {$T$-plays and $T$-strategies}
 Let us denote by $\cal T$ the set of all sequences of functions
$$T(x)=(T_1(x), T_2(x), \dots, T_k(x), \dots),$$
satisfying the following conditions:
\begin{enumerate}
\item[1.]
 The functions $T_k(x)$ ($k\in\N$) are defined and continuous 
in the half-line $[0,+\infty)$, and continuously differentiable
in $(0,+\infty)$.
\item[2.]
$0 < T_k(x)\le 1$  ($x\ge 0$, $k\in\N$).
\item[3.]
$T'_k(x)<0$, $T_{k+1}(x)<T_{k}(x)$ ($x>0$, $k\in\N$).
\item[4.]
$T_k(0)=1$ ($k\in\N$).
\end{enumerate}
 Let $T\in \cal T$.  A pair $(\alpha, n)$ is called a {\it $T$-play}
if whenever $\alpha(t)n(t)>0$  the inequality
\begin{align}
t\le T_{n(t)}(\alpha(t))
\label{Tan}
\end{align}
holds and for the action moments the inequality holds as an equality.

 Any sequence $T\in \cal T$ determines the set of all $T$-plays,
which differ from each other in who of the players uses his resource
at every action moment prescribed by the sequence $T$ (simultaneous
actions are possible).
 The $T$-plays $(\alpha_j, n_j)$ ($j=1,2$) in which the $j$-th
player begins to use his resource after his opponents's resource
has been exhausted are called the {\it simplest} $T$-plays.
 The consumption functions of the simplest $T$-plays have the form
\begin{align}
\label{baseT1}
&\alpha_1(t)=\begin{cases}
             a,\text{ for }t\in [0,T_1(a)];\\
             \text {makes the integral }
             \Ds\int\limits_{T_1(a)}^1\log(1-P_1(t))\,d\alpha_1(t)
             \text{ divergent},
            \end{cases}
\displaybreak[0]\\
\label{baseT2}
&n_1(t)=\begin{cases}
            m,&\text{for }t\in [0,T_m(a)];\\
            i,&\text{for }t\in (T_{i+1}(a),T_i(a)],\; 1\le i\le m-1;\\
            0,&\text{for }t\in (T_1(a),1],
            \end{cases}
\displaybreak[0]\\
&\alpha_2(t)=\begin{cases}
              a, &\text{for }t\in[0,T_m(a)];\\
              T^{-1}_m(t), &\text{for }t\in(T_m(a),1],
             \end{cases}
\displaybreak[0]\\
&n_2(t)= m,\text{ for }t\in [0,1]\;(\eta_1=\eta_2=\dots=\eta_m=1).
\end{align}

\begin{lm}
\label{l:Tan1}
Let $T\in \cal T$.
 The values of the payoff functions in all $T$-plays of the duel
$G_{am}(P, A)$ coincide if and only if the following equations hold:
\begin{align}
\label{eqT1}
&\exp\lt(\int_0^x\log(1-P_1(T_k(\a)))\,d\a\rt)+\prod_{i=1}^k 
\lt(1-P_2(T_i(x))\rt)=1\\
&\text{ for all\/ }0<x\le a,\; 1\le k\le m.
\nonumber
\end{align}
In this case the common value of the payoff function in all $T$-plays
of the game is equal to
\begin{align}
\label{vam1}
v_m(a)&=
A_1-(A_1+A_2)\exp\lt(\int_0^a\log(1-P_1(T_m(\a)))\,d\a\rt)=\\
&=(A_1+A_2)\prod_{i=1}^m(1-P_2(T_i(a))-A_2.
\nonumber
\end{align}
\end{lm}

\begin{pf} \underline{Necessity}.
 Suppose that the values of the payoff functions in all $T$-plays of
the duel coincide for a certain $T\in \cal T$.
 Fix $x$ ($0< x \le a$), $k$ ($1\le k \le m$), and
put $t^*=T_k(x)$.
 Denote the simplest $T$-plays of the game $G_{xk}(P,A)$
by $(\overline{\alpha}_j, \overline{n}_j)$ ($j=1,2$)
and consider two $T$-plays $(\alpha_j,n_j)$ ($j=1,2$) of
the game $G_{am}(P,A)$ satisfying the following conditions:
\begin{align}
&\alpha_1(t)=\alpha_2(t);\;n_1(t)=n_2(t)
               \text{ for }t\in[0,t^*];
\label{Tan1}
\\
&\alpha_j(t^*)=x;\; n_j(t^*)=k,\;j=1,2;
\label{Tan2}
\\
&\alpha_j(t)=\overline{\alpha}_j(t); \;
n_j(t)=\overline{n}_j(t)\text{ for }t\in(t^*,1],\;j=1,2.
\label{Tan3}
\end{align}

 We denote the action moment vectors of the plays
$(\overline{\alpha}_j,\overline{n}_j)$
of the game $G_{xk}(P,A)$
by $\overline{\eta}^j$  ($j=1,2).$
 Let us compute the values of the payoff function in these plays:
\begin{align}
\nonumber
K(\overline{\alpha}_1,\overline{\eta}_1)&=-A_2
\left(1-\prod_{i=1}^k \lt(1-P_2(T_i(x))\rt)\rt)
+A_1\prod_{i=1}^k \lt(1-P_2(T_i(x))\rt)=\\
\label{Tan4}
&=(A_1+A_2)\prod_{i=1}^k \lt(1-P_2(T_i(x))\rt)-A_2;\displaybreak[0]\\
\nonumber
K(\overline{\alpha}_2,\overline{\eta}_2)&=A_1
\left(1-\exp\lt(\int_0^x\log(1-P_1(T_k(\a)))\,d\a\rt)\rt)-\\
\nonumber
&-A_2\exp\lt(\int_0^x\log(1-P_1(T_k(\a)))\,d\a\rt)=\\
\label{Tan5}
&=A_1-(A_1+A_2)\exp\lt(\int_0^x\log(1-P_1(T_k(\a)))\,d\a\rt).
\end{align}
 By the assumption of Lemma the values of the payoff function of
the game $G_{am}(P,A)$ in the plays $(\alpha_1,n_1)$ and
$(\alpha_2,n_2)$ coincide.
 Hence the conditions \eqref{Tan1}--\eqref{Tan3} imply that the values
of the payoff function of the game $G_{xk}(P,A)$ in the plays
$(\overline{\alpha}_1,\overline{n}_1)$ and
$(\overline{\alpha}_2,\overline{n}_2)$ are equal.
 Equating \eqref{Tan4} and \eqref{Tan5}, we get \eqref{eqT1}.

\underline{Sufficiency}.
 Suppose that for a given $T\in{\cal T}$ the equation \eqref{eqT1}
holds for all $0<x\le a$,\; $1\le k\le m$.
 Let $(\alpha, n)$ be an arbitrary $T$-play.
 We need to show that $K(\alpha, n)=v_m(a)$.
 Proceed by induction in the number of units in Player II's
resource.
 For $m=0$ the statement of Lemma is true, as $K=A_1$.
 Suppose that the equation holds for $n\le k-1$ and prove it for $n=k$.
 Let $\eta$ be the action moment vector of Player II in the play
$(\alpha, n)$.
 Set
$$
 \alpha_k=\alpha(\eta_k),\quad \alpha^k(t)=\min\{\alpha_k,\alpha(t)\}.
$$
 Then
\begin{align}
\label{Tan6}
&K(\alpha, \eta_1,\dots,\eta_k)=
A_1\varphi(\alpha, 0,\eta_k)-A_2P_2(\eta_k)
\left(1-\varphi(\alpha, 0,\eta_k)\rt)+\\
\nonumber
&+(1-P_2(\eta_k))\left(1-\varphi(\alpha, 0,\eta_k)\rt)
K(\alpha^k, \eta_1,\dots,\eta_{k-1}),
\end{align}
where $\varphi(\alpha, 0,\eta_k)$ is the probability of Player I
achieving success in the time interval $[0,\eta_k)$.
 By the formula \eqref{phi1} we have
\begin{align}
\label{Tan7}
\varphi(\alpha, 0,\eta_k)=
1-\exp\lt(\int_{\alpha_k}^a\log(1-P_1(T_k(\a)))\,d\a\rt).
\end{align}
 It follows from the inductive assumption that
\begin{align}
\label{Tan8}
K(\alpha^k, \eta_1,\dots,\eta_{k-1})&=
(A_1+A_2)\prod_{i=1}^{k-1}(1-P_2(T_i(\alpha_k))-A_2.
\end{align}
 Substituting \eqref{Tan7} and \eqref{Tan8} into \eqref{Tan6}, we get
\begin{align}
\label{Tan9}
&K(\alpha, \eta_1,\dots,\eta_k)=A_1+\\
\nonumber
&+(A_1+A_2)\exp\left(\int_{\alpha_k}^a\log(1-P_1(T_k(\a)))\,d\a\right)
\left(
\prod_{i=1}^k(1-P_2(T_i(\alpha_k))-1\rt).
\end{align}
 According to \eqref{eqT1} we have
\begin{align}
\label{Tan10}
1-\prod_{i=1}^k\lt(1-P_2(T_i(\alpha_k))\rt)=\exp\int_0^{\alpha_k}
\log(1-P_1(T_k(\a)))\,d\a.
\end{align}
 Taking into account \eqref{Tan10}, we finally conclude from
\eqref{Tan9} that
\begin{align*}
K(\alpha, \eta_1,\dots,\eta_k)=A_1-
(A_1+A_2)\exp\int_0^a\log(1-P_1(T_k(\a)))\,d\a=v_k(a).
\end{align*}
 So the statement of Lemma has been proven by induction.
\end{pf}

 Let $\{T_k(x)\}\in {\cal T}$.
 The players' strategies having the form
\begin{align*}
&\xi^T(t)=
\begin{cases}
              0,& t<T_n(\a),\\
              -1/T'_n(\a),& t=T_n(\a);
\end{cases}
\qquad\eta^T=T_n(\a),
\end{align*}
where $\a$, $n$ are the players' remaining resources at a moment $t$,
are called {\it $T$-strategies}.

\begin{thm}
\label{t:Tan1}
If a sequence $\{T_k(x)\}\in {\cal T}$ satisfies the relations
\eqref{eqT1} for all $0<x\le a,$\; $1\le k\le m$, then any pair
of $T$-strategies forms an equilibrium situation (saddle point)
in the game $G_{am}(P,A)$.
The value of the game is given by the formula
\eqref{vam1}.
\end{thm}

\begin{pf}
 Let $\eta$ be an arbitrary action moment vector of Player II.
 Suppose that Player I acts according to a $T$-strategy.
 His consumption function corresponding to $\eta$ has the form
\begin{align}
\label{Tan1.1}
&\alpha^T(t)=
\begin{cases}
              \alpha_{k+1},& t\in (\eta_{k+1},T_k(\a_{k+1}));\\
              \alpha^T_k(t),& t\in [T_k(\a_{k+1}),\eta_k],
                    \text { for } \eta_k>T_k(\a_{k+1}),
\end{cases}
\end{align}
where $\alpha^T_k(t)$ is the function defined in the segment
$[T_k(\alpha_{k+1}),1]$ and inverse to $T_k(x)$, extended to
the segment $[0,T_k(\alpha_{k+1})]$ as the constant~$\alpha_{k+1}$,
$\alpha_k=\alpha^T_k(\eta_k)$ ($k=1,2\dots,m$), $\eta_{m+1}=0$,
$\alpha_{m+1}=a$.

 We will show that if $T_k(x)$ satisfies \eqref{eqT1} for all
$0<x\le a,$\; $1\le k\le m$ then the inequality
$K(\alpha^T;\eta)\ge v_m(a)$ holds.

 First let us notice that if Player I uses a $T$-strategy then
the inequality $\eta_k>T_k(\alpha_k)$ is impossible, because starting
from the moment $T_k(\alpha_{k+1})$ Player~I spends his resource
according to the function $\alpha^T_k(t)$ making the identity
$t=T_k(\alpha^T(t))$ hold.
 In the result the next action moment of Player~II
prescribed by a $T$-strategy
is being postponed.

 If for all $1\le k\le m$ one has $\eta_k=T_k(\alpha_k)$ then
we are dealing with the simplest $T$-play \eqref{baseT1},
and according to Lemma~\ref{l:Tan1} the equation
$K(\alpha^T;\eta)= v_m(a)$ holds.

 If Player II spends his resouce before the next action moment
$T_k(\alpha^T(t))$ comes, then there exist two integers $k$ and $l$
($1\le l\le k\le m$) such that
\begin{align}
\label{Tan1.1a}
&\eta_i=T_i(\alpha_i) \text{ for } k+1\le i\le m;\quad
\eta_{l-1}=T_{l-1}(\alpha_{l-1});\\
\label{Tan1.1b}
&\eta_i<T_i(\alpha_i) \text{ for } l\le i\le k.
\end{align}
 In this case, by the definition of a $T$-strategy, Player I's
resource is not being consumed in the interval
$(\eta_{k+1},T_l(\a_{k+1}))$, that is $\alpha^T(t)=\alpha_{k+1}$.
 Define the vector $\eta^1$ as follows:
\begin{align*}
&\eta^1_i=
\begin{cases}
                   \eta_i,& i=1,2,\dots,l-1,k+1,k+2,\dots,m;\\
                   T_i(\a_{k+1}),& i=l,l+1,\dots,k.
\end{cases}
\end{align*}
 Let us compute $K(\alpha^T;\eta)$ and $K(\alpha^T;\eta^1)$ by
presenting the payoff function as the sum of three summands
corresponding to the intervals $[0,\eta_k)$, $[\eta_k, T_l(\a_{k+1}))$,
$[T_l(\a_{k+1}),1]$.
 We get:
\begin{align}
\label{Tan1.2}
&K(\alpha^T,\eta)=K_{[0,\eta_k)}-A_2\Psi(\eta_k)+
\Psi(\eta_k)\prod_{i=l}^k q(\eta_i)(A_2+K_{[T_l,1]});\\
\label{Tan1.3}
&K(\alpha^T,\eta^1)=K_{[0,\eta_k)}-A_2\Psi(\eta_k)+
\Psi(\eta_k)\prod_{i=l}^k q(T_i(\alpha_{k+1}))(A_2+K_{[T_l,1]}),
\end{align}
where $K_{[0,\eta_k)}$ and $K_{[T_l,1]}$ are the mathematical
expectations of Player I's profit in the intervals $[0,\eta_k)$
and $[T_l(\a_{k+1}),1]$ when Player I's consumption function is
$\alpha^T$ and Player II's action moment vector is $\eta$,
while $\Psi(\eta_k)$ is the probability that for these consumption
function of Player I and action moment vector of Player II both
players did not achieve success up to the moment $\eta_k$.
 Since the function $q(t)$ decreases, comparing
\eqref{Tan1.2} with \eqref{Tan1.3} and taking into account
\eqref{Tan1.1a}, \eqref{Tan1.1b} we get the inequality
$$K(\alpha^T,\eta)\ge K(\alpha^T,\eta^1).$$
 Repeating the described procedure, we construct $r$ vectors
$\eta^1,\eta^2,\dots,\eta^r$ ($r<m$) such that
\begin{align*}
&K(\alpha^T,\eta)\ge K(\alpha^T,\eta^1)\ge K(\alpha^T,\eta^2)
\ge\dots\ge K(\alpha^T,\eta^r)\\
&\text {and }
\eta^r_i=T_i(\alpha_i) \text{ for all } 1\le i\le m.
\end{align*}
 According to Lemma~\ref{l:Tan1}, we have $K(\alpha^T;\eta^r) = v_m(a)$,
and consequently $K(\alpha^T;\eta)\ge v_m(a)$.

 Now let $\alpha(t)$ be an arbitrary consumption function of Player I
and $\eta^T$ be the realization of Player II's action moment vector
corresponding to $\alpha(t)$ for a $T$-strategy of Player~II.
 Let us show that $K(\alpha;\eta^T)\le v_m(a)$.
 Denote the realization of Player I's consumption function corresponding
to $\eta^T$ for a $T$-strategy of Player I by $\alpha^T$ (it is given
by the formula \eqref{Tan1.1}).
 If for all $t\in[0,1]$ such that $\alpha(t)n(t)>0$ the function
$\alpha(t)$ coincides with $\alpha^T(t)$, then we are dealing with
a $T$-play, and by Lemma~\ref{l:Tan1} we have $K(\alpha;\eta^T)= v_m(a)$.
 Otherwise there exists $t_*\in[0,1]$ for which
$\alpha(t_*)=\alpha^T(t_*)$
and there exists $\eps>0$ such that for all $t\in(t_*,t_*+\eps)$
the inequality $\alpha(t)<\alpha^T(t)$ holds.
 The inverse inequality is impossible because Player II uses
a $T$-strategy.
 Set
$$
t^*=\sup \{t: \alpha(t')<\alpha^T(t') \text{ for all } t'\in(t_*,t)\}.
$$
 The segment $[0,1]$ contains at most a countable set of segments
of the form $[t_*,t^*]$.
 Let us enumerate such segments, and let $[t_1, t'_1]$ be the first
of them.
 Define $\alpha_1(t)$ as follows:
\begin{align*}
&\alpha_1(t)=
\begin{cases}
             \alpha^T(t), &\text {for } t\in(t_1, t'_1);\\
             \alpha(t),  &\text {for } t\notin(t_1, t'_1).
\end{cases}
\end{align*}
 Let $t_1\in[\eta^T_k,\eta^T_{k-1}]$; then by the definition of
Player I's $T$-strategy we have $t'_1\in[\eta^T_k,\eta^T_{k-1}]$.
 Let us compute $K(\alpha;\eta^T)$ and $K(\alpha_1;\eta^T)$ by
presenting the payoff function as the sum of three summands
corresponding to the intervals $[0,t_1)$, $[t_1, t'_1)$, $[t'_1,1]$.
 We get:
\begin{align}
\label{Tan1.4}
K(\alpha,\eta^T)&=K_{[0,t_1)}+A_1\Psi(t_1)
-\Psi(t_1)(A_1-K_{[t_1',1]})\exp\int
\limits_{t_1}^{t'_1}\mu(t)\,d\alpha(t);\\
\label{Tan1.5}
K(\alpha_1,\eta^T)&=K_{[0,t_1)}+A_1\Psi(t_1)
-\Psi(t_1)(A_1-K_{[t_1',1]})\exp\int
\limits_{t_1}^{t'_1}\mu(t)\,d\alpha^T(t),
\end{align}
where $K_{[0,t_1)}$ and $K_{[t_1',1]}$ are the mathematical expectations
of Player I's profit in the intervals $[0, t_1)$ and $[t_1',1]$ when
Player I's consumption function is $\alpha(t)$ and Player II's action
moment vector is $\eta^T$, while $\Psi(t_1)$ is the probability that for
these consumption function of Player I and action moment vector of
Player II both players did not achieve success up to the moment
$t_1$, and $\mu(t)=-\log p(t)$.
 Note that
\begin{align}
\label{Tan1.6}
\int\limits_{t_1}^{t'_1}\mu(t)\,d\alpha(t)>
\int\limits_{t_1}^{t'_1}\mu(t)\,d\alpha^T(t).
\end{align}
Indeed, integrating by parts we get
\begin{align}
\label{Tan1.7}
&\int_{t_1}^{t'_1}\mu(t)\,d\alpha(t)
=\mu(t'_1)\int_{t_1}^{t'_1}\,d\alpha(t)-
\int_{t_1}^{t'_1}
\left(\int_{t_1}^{t}\,d\alpha(\tau)\right)\,d\mu_1(t),\\
\label{Tan1.8}
&\int_{t_1}^{t'_1}\mu(t)\,d\alpha^T(t)
=\mu(t'_1)\int_{t_1}^{t'_1}\,d\alpha^T(t)-
\int_{t_1}^{t'_1}
\left(\int_{t_1}^{t}\,d\alpha^T(\tau)\right)\,d\mu_1(t).
\end{align}
By the definition of the segment $[t_1,t'_1]$, for all $t\in[t_1,t'_1]$
the following inequality holds:
\begin{align}
\label{Tan1.9}
           \int_{t_1}^{t}\,d\alpha(\tau)\le
           \int_{t_1}^{t}\,d\alpha^T(\tau);
\end{align}
moreover, the inequality turns into an equality for $t=t'_1$ only.
 Using \eqref{Tan1.9}, one deduces \eqref{Tan1.6} from \eqref{Tan1.7}
and \eqref{Tan1.8}.
 Comparing $K(\alpha; \eta^T)$ with $K(\alpha_1; \eta^T)$
(the formulas \eqref{Tan1.4}, \eqref{Tan1.5}) and taking
\eqref{Tan1.6} into account, we get:
\begin{align*}
 K(\alpha,\eta^T)\le K(\alpha_1,\eta^T).
\end{align*}

 Repeating the described procedure, we construct a sequence of
functions $\alpha_k$ such that
\begin{align*}
 K(\alpha_k,\eta^T)\le K(\alpha_{k+1},\eta^T),\;k\in\N.
\end{align*}
 So
\begin{align}
\label{Tan1.10}
 K(\alpha,\eta^T)\le K(\alpha_{k},\eta^T) \text{ for any }k\in\N.
\end{align}
 Let $\alpha^*(t)=\Ds\lim_{k\to\infty}\alpha_k(t)$.
 Passing to the limit for $k\to\infty$ in the inequality \eqref{Tan1.10}
and using Helly's convergence theorem~\cite{kf} we conclude that
\begin{align*}
 K(\alpha,\eta^T)\le K(\alpha^*,\eta^T).
\end{align*}
 Since $\alpha^*(t)=\alpha^T(t)$ for all $t\in\{t:n(t)>0\}$,
by Lemma~\ref{l:Tan1} we have $K(\alpha^*;\eta^T)= v_m(a)$ and thus
$K(\alpha;\eta^T)\le v_m(a)$.
\end{pf}

\begin{cor}
\label{c:Tan1}
 If the function $P_2(t)$ strictly increases in the segment $[0,1]$,
then there exists at most one sequence $\{T_k(x)\}\in {\cal T}$
satisfying \eqref{eqT1} for all $x\ge 0$, $k\in\N$.
\end{cor}

\begin{pf}
 Suppose there exist two sequences $T^1,\;T^2\in{\cal T}$
satisfying \eqref{eqT1} for all $x\ge 0$, $k\in\N$.
 Let
\begin{align}
\label{cor.1}
l=\min\{k:T_k^1\neq T_k^2\} \quad\text{and}\quad
T_l^1(a)\neq T_l^2(a),\; a>0.
\end{align}
 By Theorem~\ref{t:Tan1} the game $G_{al}(P,A)$ has the value equal
to the value of the payoff function in the $T^1$- and $T^2$-plays,
that is the following equation holds:
\begin{align}
\label{cor.2}
(A_1+A_2)\prod_{i=1}^{l}\lt(1-P_2(T^1_i(a))\rt)-A_2=
(A_1+A_2)\prod_{i=1}^{l}\lt(1-P_2(T^2_i(a))\rt)-A_2.
\end{align}
 But by the definition of $l$ for all $i<l$ one has $T_i^1(a)=T_i^2(a)$,
hence using \eqref{cor.2} and taking into account the strict
monotonicity of the function $P_2(t)$ we conclude that
$T_l^1(a)=T_l^2(a)$.
 We have come to a contradiction which proves uniqueness
of the sequence $\{T_k(x)\}$.
\end{pf}

 Let $\{T_k(x)\}\in {\cal T}$ be a sequence satisfying the relation
\eqref{eqT1} for all $a\ge 0$ and $k\in \N$.
 Introduce the notation
\begin{align}
\label{eqT2.1a}
\pi_0(x)=1;\;
 \pi_{k}(x)=q(T_{k}(x))\pi_{k-1}(x)\;(k\in\N).
\end{align}
 Then
\begin{align}
\exp\lt(\int_0^x\log(p(T_k(\a))\,d\a\rt)=1-\pi_{k}(x).
\label{eqT2.1}
\end{align}
 Differentiating \eqref{eqT2.1} in $x$, we get
\begin{align}
\label{eqT2.2}
\log(p(T_k(x))=-\frac{\pi'_{k}(x)}{1-\pi_{k}(x)}
\end{align}
 Let us write down the recurrence relation for $\pi'_{k}(x)$:
\begin{align}
\label{eqT2.3}
\pi'_{k}(x)=q'(T_k(x))T'_k(x)\pi_{k-1}(x)+q(T_k(x))\pi'_{k-1}(x).
\end{align}
 It follows from \eqref{eqT2.2} and \eqref{eqT2.3} that the sequence
$\{T_k(x)\}$ satisfies the system of ordinary differential equations
\begin{align}
\label{difeq}
&\frac{dT_k}{dx}=\phi(T_1,T_2, \dots, T_k), \quad
x\ge 0;\;  k\in \N,
\end{align}
where
\begin{align*}
&\phi(T_1,T_2, \dots, T_k)=\\
&=-\frac{\left(
1-\Ds\prod_{i=1}^k q(T_i)\rt)\log p(T_k)-q(T_k)
\lt(1-\Ds\prod_{i=1}^{k-1}q(T_i)\rt)\log p(T_{k-1})}
{q'(T_k)\Ds\prod_{i=1}^{k-1}q(T_i)}.
\end{align*}

\begin{lm}
\label{l:Tan2}
 Assume that $p(t)$ and $q(t)$ are continuously differentiable in
$(0,1]$, $p(0)=q(0)=1$; $p(1)=q(1)=0$; $p(t)>0$ for $t<1$;
$p'(t)\le 0$; $q'(t)<0$.
 Then the system of ordinary differential equations \eqref{difeq}
under the initial conditions
\begin{align}
T_k(0)=1 \quad (k\in N)
\label{difeq1}
\end{align}
has a solution in the half-line $x>0$; moreover,
$\{T_k(x)\}\in {\cal T}$.
\end{lm}

\begin{pf}
 Let us prove Lemma by induction in the number of action moments
of Player II.
 The first equation of the system \eqref{difeq} has the form:
\begin{align}
\label{difeq2}
&\frac{dT_1}{dx}=
-\frac{(1-q(T_1))\log p(T_1)}{q'(T_1)}.
\end{align}
 Integrating \eqref{difeq2} under the initial condition $T_1(0)=1$,
we get:
\begin{align}
\label{difeq3}
x(T_1)=
\int\limits_{T_1}^1\frac{q'(\tau)\,d\tau}{(1-q(\tau))\log p(\tau)}.
\end{align}
 The function $T_1(x)$ is the inverse function to $x(T_1)$.
 Let us check that it satisfies the conditions 1--4.

 First we have to show that $T_1(x)$ is defined in the half-line
$[0,+\infty)$.
 Choose $\delta>0$ such that $-\log p(\delta)<1$. Then
$$
x(t)\ge
\int\limits_{t}^{\delta}\frac{(1-q(\tau))'\,d\tau}{(1-q(\tau))}+
x(\delta) = \log(1-q(\delta))-\log(1-q(t))+x(\delta).
$$
 Since the right hand side of the inequality tends to $+\infty$ as
$t\to +0$, we have $T_1\to +0$ as $x\to +\infty$ and therefore
the function $T_1(x)$ is defined for all $x>0$.
 According to the initial condition, $T_1(0)=1$.

 By \eqref{difeq2}, one has $T'_1(x)<0$.
 So $T_1(x)$ decreases from 1 for $x=0$ to $0$ as $x\to +\infty$.
 Hence $0 < T_1(x)\le 1$.

 Suppose that for $1\le i\le k-1$ a solution sequence $T_i(x)$ of
the system \eqref{difeq}, \eqref{difeq1} exists and satisfies
the conditions 1--4.
 Substituting it into the $k$-th equation of the system
\eqref{difeq}, we get:
\begin{align}
\label{difeq4}
\frac{dT_k}{dx}=\Phi_k(T_k,x),
\end{align}
where
\begin{align}
\label{difeq5}
\Phi_k(t,x)=\phi_k(T_1(x),T_2(x),\dots, T_{k-1}(x),t).
\end{align}
 Let us show that the equation \eqref{difeq4} under the initial
condition
\begin{align}
\label{difeq6}
T_k(0)=1
\end{align}
has a solution $T_k(x)$, which satisfies the monotonicity
condition
\begin{align}
\label{difeq7}
T'_k(x)<0,\; T_{k}(x)<T_{k-1}(x) \text{ for all }x>0.
\end{align}
 Note that the following inequalities holds for all $x>0$: 
\begin{align}
\label{difeq8}
\Phi_k(T_{k-1}(x),x)\le\frac{T'_{k-1}(x)}{q(T_{k-1}(x))}<T'_{k-1}(x)<0.
\end{align}
 Indeed,
\begin{align}
\label{difeq9}
\Phi_k(T_{k-1}(x),x)=
-\frac{\left(1-q(T_{k-1}(x))\rt)\log p(T_{k-1}(x))}
{q'(T_{k-1}(x))\pi_{k-1}(x)},
\end{align}
where $\pi_{k}(x)$ are the functions defined by
the formulas \eqref{eqT2.1a}.
 On the other hand, by the inductive assumption
$$T_{k-1}(x)<T_{k-2}(x),$$
hence \eqref{difeq} implies the inequality
\begin{align}
\label{difeq10}
T'_{k-1}(x)\ge
-\frac{\left(1-q(T_{k-1}(x))\rt)\log p(T_{k-1}(x))}
{q'(T_{k-1}(x))\pi_{k-2}(x)}.
\end{align}
 Comparing \eqref{difeq9} and \eqref{difeq10}, we obtain \eqref{difeq8}.
 Denote the  numerator of the fraction in the right hand side of
the equation \eqref{difeq} by $F_{k}(t,x)$, that is
\begin{align*}
F_{k}(t,x)= (1-\pi_{k-1}(x)q(t))\log p(t)-
(1-\pi_{k-1}(x))q(t)\log p(T_{k-1}(x)).
\end{align*}
 Next we will show that the equation
\begin{align}
\label{difeq11}
F_{k}(t,x)=0,\;x>0,
\end{align}
determines an implicit function $t=f_k(x)$, which has the following
properties:
\begin{enumerate}
\item
the function $f_k(x)$ ($k\in\N$) is defined and continuously
differentiable in the half-line $(0,+\infty)$;
\item
$f_k(x)\to 1$ as $x\to +0$;
\item
$f'_k(x)<0$\;  ($x>0$);
\item
$f_k(x)<T_{k-1}(x)$\; ($x>0$, $k\ge 2$).
\end{enumerate}
 It follows from \eqref{difeq8} that
$$F_k(T_{k-1}(x),x)<0 \text{ for all } x>0.$$
 On the other hand, for any $x>0$
$$
\lim_{t\to +0}F_{k}(t,x)=
-(1-\pi_{k-1}(x))\log p(T_{k-1}(x))>0.
$$
 Therefore, for any $x>0$ the equation \eqref{difeq11} has a solution
$$t=f_k(x)\in(0,T_{k-1}(x)).$$
 To prove that the solution is unique let us check that
$$\frac{\partial F_{k}}{\partial t}<0 
\text { for all } x>0,\;t\in(0,1).$$
 Indeed,
\begin{align*}
\frac{\partial F_{k}}{\partial t}&=
\frac{p'(t)}{p(t)}
(1-\pi_{k-1}(x)q(t))-
q'(t)\pi_{k-1}(x)\log p(t)
-\\&
-q'(t)\pi_{k-1}(x)\log p(T_{k-1}(x))<0.
\end{align*}
 It follows from the implicit function theorem that the equation
\eqref{difeq11} determines an implicit function $t=f_{k}(x)$,
which is differentiable in the half-line $x>0$.
 Let us check that $f'_{k}(x)<0$.
 We have shown that $F'_{k}(t)<0$, so in view of the relation
$$ (f_k)'_x=-\frac{(F_{k})'_t}{(F_{k})'_x}$$
it suffices to check that
$$ (F_{k})'_x<0 \text{ for all } x>0,\; t<T_{k-1}(x).$$
 Taking into account the fact that, according to \eqref{eqT2.2},
$$\pi'_{k-1}(x)=-(1-\pi_{k-1}(x))\log(p(T_{k-1}(x)),
$$
we have
\begin{align*}
 (F_{k})'_x&=
-(1-\pi_{k-1}(x))\log(p(T_{k-1}(x))
\left(\log p(t)-\log(p(T_{k-1}(x))\rt)-\\
&-\frac{p'(T_{k-1}(x))}{p(T_{k-1}(x))}T'_{k-1}(x)(1-\pi_{k-1}(x))q(t).
\end{align*}
 Since $t<T_{k-1}(x)$ and $T'_{k-1}(x)<0$, it follows that
$(F_{k})'_x<0$.

 It remains to check that $f_k(x)\to 1$ as $x\to +0$.
 It was proven above that $f_k(x)$ decreases monotonically
in the half-line $x>0$.
 Taking into account the inequalities
$$f_k(x)<T_{k-1}(x)<1,$$
we conclude that there exists a limit of $f_k(x)$ as $x\to +0$ and
$$\lim_{x\to +0}f_k(x)=c\le 1.$$
 Suppose that $c<1$.
 Substituting $t=f_k(x)$ into \eqref{difeq11}, we get:
\begin{align}
\label{difeq12}
(1-\pi_{k-1}(x)q(f_k(x)))\log p(f_k(x))=
(1-\pi_{k-1}(x))q(f_k(x))\log p(T_{k-1}(x)),
\end{align}
 As $x\to +0$, the right hand side of the equation \eqref{difeq12}
tends to $-\infty$, and the limit of the left hand side is equal to
$\log p(c)>-\infty$.
 It follows from this contradiction that
$$\lim_{x\to +0}f_k(x)=1.$$

 Now we are ready to proceed with the construction of the function
$T_k(x)$.
 Associate with any $a>0$ the solutions $y_a(x)$ and $z_a(x)$
of the equation \eqref{difeq4} in the half-line $x\ge a$
satisfying the initial conditions
$$y_a(a)=T_{k-1}(a);\quad z_a(a)=f_{k}(a).
$$
 This equation in the domain $x>0$, $0<t<1$ satisfies the conditions
of the theorem on the existence and uniqueness of solutions.
 By the inequalities \eqref{difeq8} and $f'_{k-1}(x)<0$, the curves
$y_a(x)$ and $z_a(x)$ for $x>a$ are situated between the curves
$$t=f_{k}(x) \text{ and } t=T_{k-1}(x).
$$
 Take $c>0$ and denote $\Ds\inf_{a>0}y_a(c)$ by $b$.
 Let $\tilde{y}(x)$ be a solution of the equation \eqref{difeq4}
satisfying the initial condition $\tilde{y}(c)=b$.
 By the uniqueness theorem, for any $a>0$ the integral curve
$\tilde{y}(x)$ is situated strictly between curves $y_a(x)$
and $z_a(x)$ for all $x>0$.
 Hence $\tilde{y}(x)$ can be extended to the half-line $x\ge 0$;
moreover, $\tilde{y}(0)=1$ and
$$f_{k}(x)<\tilde{y}(x)<T_{k-1}(x) \text{ for all } x>0.
$$
 Since $F_k(t,x)$ decreases in $t$ for any $x>0$, one has
$$F_k(\tilde{y}(x),x)<F_k(f_{k}(x),x)=0,$$
and therefore $\tilde{y}'(x)<0$.
 Thus $\tilde{y}(x)$ satisfies all the conditions imposed on
the function $T_k(x)$, and so the existence of this function
is proven. Put $T_k(x)=\tilde{y}(x)$.
\end{pf}

\begin{lm}
\label{l:Tan2a}
 Assume that $p(t)$ and $q(t)$ are continuously differentiable in
$(0,1]$, $p(0)=q(0)=1$; $p(1)=q(1)=1$; $p(t)>0$ for $t<1$;
$p'(t)\le 0$; $q'(t)<0$.
 Then a solution of the system of differential equations
\eqref{difeq}, \eqref{difeq1} is unique.
\end{lm}

\begin{pf}
 Let us proceed by induction.
 For $k=1$ the function $T_1(x)$ is inverse to the function $x(T_1)$,
which is determined uniquely by the formula \eqref{difeq3}.
 Suppose that for $1\le i\le k-1$  the system \eqref{difeq},
\eqref{difeq1} has a unique solution $T_i(x)$.
 Let us show that the problem \eqref{difeq4}, \eqref{difeq6}
has a unique solution.
 Let
$T_k(x)$ be the solution of this problem constructed
in Lemma~\ref{l:Tan2} and
$y(x)$ be an arbitrary solution of this problem.
 We will show that
$$y(x)\equiv T_k(x).$$
 Consider two cases.
\begin{enumerate}
\item[1.]
 The integral curve $y(x)$ of the equation \eqref{difeq4} for all $x>0$
is situated strictly between the graphs of the functions
$$y=f_{k}(x) \text{ and } y=T_{k-1}(x).
$$
 Then, according to the above, $y(x)$ satisfies the conditions
\eqref{difeq6}, \eqref{difeq7}, and by Corollary~\ref{c:Tan1},
$$y(x)\equiv T_k(x).
$$
\item[2.]
For some $x>0$ one of the inequalities
$$f_{k}(x)<y(x)<T_{k-1}(x)
$$
is false.
 We will show that in this case the curve $y(x)$ does not go through
the point $(0,1)$, i.~e., in this case the initial condition
\eqref{difeq6} is not satisfied.
\begin{enumerate}
\item
 Suppose there exists $x_0>0$ such that $y(x_0)<f_{k}(x_0)$.
Then $y'(x)>0$ for all $x\in(0,x_0)$, and therefore
$y(x)<f_{k}(x_0)$ for all $0<x<x_0$, so in particular $y(0)<1$.
\item
 Suppose there exists $x_1>0$ such that $y(x_1)>T_{k-1}(x_1)$.
 Then by the inequality \eqref{difeq8} one has
$y(x)>T_{k-1}(x)$ for all $0<x<x_1$.
 Since $q'(t)$ is continuous and $q'(t)<0$ for $t\in(0,1]$,
there exist two numbers $c_1$ and $c_2$ such that
$$c_1\le -q'(t)\le c_2 \text { for } t\in [T_{k-1}(x_1),1].
$$
 It follows from the function $F_k(t,x)$ being monotonically decreasing
in $t$ together with the ineguality \eqref{difeq8} that for
$t>T_{k-1}(x)$ one has
\begin{align}
\label{difeq13}
\Phi_k(t,x)=-\frac{F_k(t,x)}{q'(t)\pi_{k-1}(x)}<
-\frac{F_k(T_{k-1}(x),x)}{q'(t)\pi_{k-1}(x)}<
\frac{T'_{k-1}(x)q'(T_{k-1}(x))}{q'(t)q(T_{k-1}(x))}.
\end{align}
 We find $\delta>0$ such that $q(T_{k-1}(\delta))<c_1/c_2$ and put
$\delta_1=\min(x_1,\delta)$.
 Then for $x\in(0,\delta)$, $t>T_{k-1}(x)$ the inequality
$$\Phi_k(t,x)<T'_{k-1}(x)$$
holds, hence there exists $\varepsilon>0$ such that for $x\in(0,\delta_1)$
the inequality $y(x)-T_{k-1}(x)>\varepsilon$ is satisfied.
 Thus $y(x)$ does not go through the point $(0,1)$.
\end{enumerate}
\end{enumerate}
\end{pf}

\begin{rem}
\label{Rem1}
 Suppose $P_2(t)=t$; and let $y(x)$, $z(x)$ be the solutions of
the equation \eqref{difeq4} in the half-line $[a,+\infty)$ under
the initial conditions $y(a)=y_0$,  $z(a)=z_0$,
where
$$f_k(a)\le z_0<y_0\le T_{k-1}(a).$$
 Then the difference $e(x)=y(x)-z(x)$ decreases in $x$.
\end{rem}

\begin{pf}
 Consider the derivative of the difference:
$$
e'(x)=\Phi_k(y(x),x)-\Phi_k(z(x),x).
$$
 From the relation
$$
\left(\Phi_k(t,x)\right)'_t=\left(F_k(t,x)\right)'_t/\pi_{k-1}(x),
$$
taking into account the inequality $(F_k)'_t<0$ obtained in the proof
of Lemma~\ref{l:Tan2}, we get $\left(\Phi_k\right)'_t<0$ for all $x>0$.
 Hence $\Phi_k(t,x)$ decreases in $t$ for any $x>0$.
 Since $y(x)>z(x)$, we have $e'(x)<0$ and therefore the difference
$y(x)-z(x)$ decreases in $x$.
\end{pf}

 The next Theorem follows from Theorem~\ref{t:Tan1} and
Lemma~\ref{l:Tan2}.

\begin{thm}
\label{t:Tan2}
 If in the noisy fighter-bomber duel the players' accuracy functions
$P_j(t)$ ($j=1,2$) are continuously differentiable in the segment
$(0,1]$, $P_j(0)=0$, $P_j(1)=1$, $P_j(t)<1$ for $t<1$, $P'_1(t)\ge 0$,
$P'_2(t)>0$, then the game has an equilibrium situation (saddle point)
in pure strategies.
 The value and the optimal strategies have the form:
\begin{align}
\label{fin}
&v_m(a)=(A_1+A_2)\prod_{i=1}^m(1-P_2(T_i(a))-A_2;\\
&\xi^T(t)=
\begin{cases}
              0,& t<T_n(\a),\\
              -1/T'_n(\a),& t=T_n(\a);
\end{cases}
\quad\eta^T=T_n(\a),
\end{align}
where $\a$, $n$ are the players' remaining resources at the moment $t$,
and $\{T_k(x)\}$ is a solution of the system of equations
\eqref{difeq}, \eqref{difeq1}.
\end{thm}

\section{Conclusions}

\begin{enumerate}
\item
 The optimal strategies $\xi^T$, $\eta^T$ prescribe the players
to refrain from using their resources until the moment $T_m(a)$
(let us call it the good start-consume moment).
\item
 The good start-consume moment is a function  $T(\alpha,\beta)$ of
the current values $\alpha$ and $\beta$ of the player's resources;
it does not depend on the development of the game up to the current
moment and is common for both players.
\item
 The more resources the players have, the earlier the good
start-consume moment comes, because $T_k(x)$ decreases in $x$ and $k$.
\item
 The optimally behaving players begin using their resources
at the good start-consume moment (one of the players or both of them).
\item
 If at the good start-consume moment Player I starts to act according
to the optimal consumption function, then he continues to act until
the next Player II's action moment and during the whole of this period
of time the following equation holds:
\begin{equation}
t=T_{m(t)}(\a(t)).
\label{Concl}
\end{equation}
\item
 If at the good start-consume moment Player II acts, then his resource
decreases by one unit and the players stop using their resources until
the good start-consume moment corresponding to the current values
of resources.
\item
 If Player I behaves optimally, the optimal strategy of Player II
prescribes him to act at any moment when the equation \eqref{Concl}
holds (possibly simultaneously with Player I, interrupting his
actions), or to refrain from using his resource until the end of
the game (moment $t=1$), which does not affect the payoff.
\item
 If one of players uses a $T$-strategy, then a play is realized
satisfying the condition
$t\le T_{m(t)}(\a(t))$  for all $t\in[0,1]$
such that $\a(t)m(t)>0$.
\item
 If both players act according to $T$-strategies, then they consume
their resources only at those moments of time for which the equation
\eqref{Concl} holds.
 In this case one of the $T$-plays is realized.
 The payoff function takes the same values at all of these plays.
\end{enumerate}

 Note that the sequence of functions $T_k(x)$ is a continuous analogue
of the infinite matrix $\{t_{mn}\}$ ($m\in\N$, $n\in\N$) of
``good first-shot times'' of the noisy duel with discrete resources of
both players \cite{fk}.

\section{Appendix. Numerical solution of the game}

 We will consider the case when $P_2(t)=t$, which does not restrict
the generality.
 Indeed, let $P_2(t)\not\equiv t$, $P'_2(t)>0$ and $P_2(t)>0$ for
$t\in (0,1]$.
 Let us make the change of variables $\tau=P_2(t)$ and solve the game
$G_{am}(\tilde P,A)$, where $\tilde P=(P_1(P_2^{-1}(\tau)),\tau)$.
 Obviously the values of the games $G_{am}(P,A)$ and
$G_{am}(\tilde P,A)$ are equal and the optimal $T$-strategies of
the game $G_{am}(P,A)$ are determined by the sequence
$T_k(x)=P_2^{-1}(\tilde T_k(x))$, where $\tilde T_k(x)$ is
the solution of the problem \eqref{difeq}, \eqref{difeq1} for
the game $G_{am}(\tilde P,A)$.

 By \eqref{difeq2} $T'_1(x)\to -\infty$ as $x\to +0$, so it follows from
the inequality \eqref{difeq7} and the initial conditions \eqref{difeq1}
that $\Ds\liminf_{x\to +0}T'_k(x)=-\infty$ for all $k\in \N$.
 Hence the system of equations \eqref{difeq} has a singularity at
the point $x=0$, and therefore it is impossible to solve the Cauchy
problem for this system with the initial conditions at the point $x=0$.
 We will integrate the system \eqref{difeq} using the method described
in the proof of Lemma~\ref{l:Tan2}, that is we will find solutions
$y_k(x)$ and $z_k(x)$ of the equation \eqref{difeq4} in the segment
$[\delta_k,a] $ under the initial conditions
\begin{equation}
y_k(\delta_k)=T_{k-1}(\delta_k),\quad
z_k(\delta_k)=f_k(\delta_k),
\end{equation}
where $\delta_k>0$ is a small number and $f_k(x)$ is the implicit
function determined by the equation \eqref{difeq11}.
 The curves $y_k(x)$ and $z_k(x)$ (we call will them the $k$-th upper
and the $k$-th lower solutions) bound the desired curve from above
and from below:
$$
  z_k(x)<T_k(x)<y_k(x).
$$
 By Remark~\ref{Rem1}, the difference $y_k(x)-z_k(x)$ decreases in $x$
and therefore the following estimate holds:
\begin{equation}
\label{deltak}
\Delta_k= \max_{x\ge\delta_k} |y_k(x)-z_k(x)|=
|y_k(\delta_k)-z_k(\delta_k)|=
T_{k-1}(\delta_k)-f_{k}(\delta_k).
\end{equation}
 From the continuity of the functions $T_{k-1}(x)$, $f_k(x)$ and
the equations $T_{k-1}(0)=f_k(0)=1$ it follows that
$T_{k-1}(\delta_k)-f_k(\delta_k)\to 0$ as $\delta_k\to 0$.

 To find the function $T_1(x)$ one needs to tabulate the function
$$
x(T_1)=
-\int\limits_{T_1}^1\frac{dt}{t\log(1- P_1(t))}.
$$
and find the inverse function.
 Using the tabulated approximate values of the function $T_1(x)$
in the subsequent computations is undesirable, since in
the computation of $T_1(x)$ in a neiborhood of the point $x=0$
we lose precision.
 The right hand side of the system \eqref{difeq} does not depend
on $x$ explicitly.
 It depends on $T_k(x)$ only, so the change of variables $u=T_1(x)$
allows to solve the further equations of the system ($k\ge 2$)
without using $T_1(x)$.
 Set
$$
\tilde T_k(u)= T_{k+1}\left(T_1^{-1}(u)\right),
 \quad k=1,2\dots, m-1,
$$
then
$$
\left.\tilde T'_k(u)= T'_{k+1}(x)/T'_1(x)\right|_{x=T_1^{-1}(u)}.
$$
 Since $T'_1(x)\to -\infty$ as $x\to +0$, passing to the variable
$u=T_1(x)$ decreases the absolute values of derivatives of
the functions we are looking for, which increases the precision of
the computations.
 After the change of variables we get a system of differential
equations in the segment $[T_1(a),1]$:
\begin{equation*}
\frac{d\tilde T_k}{du}=
\frac{\tilde \phi_k(u,\tilde T_1,\tilde T_2,\dots,
\tilde T_k)}{u\log p(u)}
\end{equation*}
under the initial conditions
$$
 \tilde T_k(1)=1,  \quad k=1,2\dots m-1.
$$
 Suppose that the first $k-1$ functions
$$
 \tilde T_i(u),  \quad i=1,2\dots, k-1
$$
have been found.
 Then $\tilde T_k(u)$ is the solution of the problem
\begin{align}
\label{difeqp1}
&\frac{d\tilde T_k}{du}=
\tilde \Phi_k(\tilde T_k,u)\\
& \tilde \Phi_k(t,u)=
\frac{\tilde \phi_k(u,\tilde T_1,\tilde T_2,\dots,\tilde T_{k-1},t)}
{u\log p(u)}
\end{align}
under the initial condition
$$
 \tilde T_k(1)=1.
$$
 After the change of variables, the initial conditions for the upper
and lower curves take the form:
\begin{align*}
&\tilde y_k(u_k)=\tilde T_{k-1}(u_k);\quad
\tilde z_k(u_k)=\tilde f_k(u_k);\\
&\tilde T_0(u_k)=u_k;\quad k=1,2,\dots,n;
\end{align*}
where $u_k=1-\delta_k$, $\delta_k>0$, and $\tilde f_k(u_k)$
is the solution of the equation $\tilde \Phi_k(t,u_k)=0$
with respect to $t$.
 In view of the strict monotonicity of the function $T_1(x)$,
it follows from \eqref{deltak} that
$$
 \max_{u\in[T_1(a),u_k]} |\tilde y_k(u)-\tilde z_k(u)|=
 \tilde T_{k-1}(u_k)-\tilde f_k(u_k)\to 0 
 \text { as } u_k\to 1-0\;(k\ge1).
$$

 Let us briefly describe a numerical algorithm for solving the game.
 The purpose of the algorithm is to compute the value of the game
$G_{ak}(P,A)$ where $P(t)=(P_1(t),t)$  $k=1,2,\dots,m$ and
tabulate the functions $T_k(x)$ in the segment $[a_0,a]$, $a_0>0$
with a given step $h$.
 The algorithm's work consists of two stages.

 Stage 1. Compute the values of the function $T_1(x)$ in the segment
$[a_0,a]$ with the step $h$ by solving the equation
$$
x(T_1)=a_0+(i-1)h;\quad i=1,2\dots, M_a,\;
M_a=\left[\frac{a-a_0}{h}\right]+1,
$$
where $x(t)$ is the function defined by the formula
\begin{align}
x(t)=-\int\limits_{t}^1\frac{d\tau}{\tau\log p(\tau)}.
\end{align}

 Stage 2. Compute the values of the function $T_k(x)$ in the segment
$[a_0,a]$ with the step $h$ ($k=2,\dots,m$).
 At the level $k$ for tabulating the function $T_k(x)$ one performs
the following computations:
\begin{enumerate}
\item
 Tabulate the $k$-th upper solution $\tilde y_k(u)$ of the equation
\eqref{difeqp1} in the segment $[T_1(a),u_0]$ ($u_0<1$) under initial
condition $\tilde y_k(u_0)=\tilde T_{k-1}(u_0)$.
\item
 Find an approximate solution of the equation $\tilde \Phi_k(t,u_0)=0$
with respect to $t$ in the segment $[0,u_0]$.
 Denote the solution of this equation by $\tilde f_k$.
\item
 Tabulate the $k$-th lower solution $\tilde z_k(u)$ of the equation
\eqref{difeqp1} in the segment $[T_1(a),u_0]$ under the initial
condition $\tilde z_k(u_0)=\tilde f_k$.
\item
 Tabulate the function $\tilde T_k(u)$ by the formula
$$
\tilde T_k(u)=(\tilde y_k(u)+\tilde z_k(u))/2
$$
in the segment $[T_1(a),u^*]$, where
$$u^*=\Ds\max\{u<u_0:\tilde y_k(u)-\tilde z_k(u)<\varepsilon\},\;
\varepsilon \text { is the given precision}. $$
\item
 Returning to original variable $x$, tabulate the function $T_{k}(x)$
in the segment $[a_0,a]$.
\item
 Compute the value of the game $v_k(a)$ by the formula \eqref{fin}.
\end{enumerate}

\end{document}


\newcommand{\lt}{\left }
\newcommand{\rt}{\right }
\renewcommand{\le}{\leqslant}
\renewcommand{\ge}{\geqslant}

\newcommand{\R}{\bf R }
\newcommand{\N}{\bf N}
\newcommand{\Z}{\bf Z}


\newcommand{\D}{\Delta}
\renewcommand{\d}{\delta}
\newcommand{\g}{\gamma}
\newcommand{\G}{\Gamma}
\newcommand{\Th}{\varTheta}
\renewcommand{\l}{\lambda}
\newcommand{\eps}{\varepsilon}
\newcommand{\ff}{\varphi}
\renewcommand{\a}{\alpha}
\renewcommand{\b}{\beta}

\newtheorem{thm}{Теорема}
\newtheorem{lm}[thm]{Лемма}
\newtheorem{cor}[thm]{Следствие}
\newtheorem{rem}[thm]{Замечание}
\newenvironment{pf}{\par\medskip\noindent{\bf Доказательство.}}%
{\nopagebreak\hfill$\Box$\par\medskip}

\title{ШУМНАЯ ДУЭЛЬ ПУЛЕМЕТЧИКА \\ СО СНАЙПЕРОМ}
\author{Посицельская Л.Н.}
\date{}
\maketitle

\begin{abstract}

\vskip 1.25mm
Рассматривается игра типа дуэли, в которой первый игрок расходует свой ресурс
непрерывно, а второй --- дискретными порциями. Каждому игроку известны текущие
значения ресурсов обоих игроков.
 Дано решение  игры  в явном виде.
\vskip 1,5mm

Ключевые слова: шумная дуэль, платежная функция, страте\-гия,
цена игры, расход ресурса.
\end{abstract}

\section {Введение}
   Классическая дуэль есть игра двух лиц с нулевой
суммой следующего вида. Игроки располагают определенными ресурсами и
используют их в течение заданного промежутка времени с целью достижения
успеха. Применение в момент $t$ ресурса $\gamma$ приводит к успеху с
вероятностью, зависящей только от времени $t$ (обычно предполагается, что
вероятность успеха возрастает по времени) и величины ресурса $\gamma$.
Как только один из игроков достигает цели, он получает выигрыш, равный
проигрышу  соперника, и игра прекращается.
   Различные предположения о способе использования игроками своего ресурса
и о поступлении  информации о поведении противника в ходе игры порождают
разнообразные виды дуэлей~\cite{kar,kimoverview}.
   Исследованы модели, в которых ресурсы игроков являются дискретными
(дуэли {\it снайперов\/}), бесконечно делимыми
(дуэли {\it пулеметчиков\/}), дуэли с непрерывным расходом ресурса
у одного из игроков и дискретным у другого, называемые  {\it смешанными\/}
дуэлями, или дуэлями {\it пулеметчика со снайпером\/} \cite{pdn,pnr}.
   Изучались {\it шумные\/} дуэли \cite{pdn,fk}, в которых каждый из игроков
в данный момент времени располагает информацией о поведении противника
до этого момента, и {\it бесшумные\/} дуэли, в которых не предполагается
поступление такой информации.
   В настоящее время дуэли считаются классическими моделями конкурентной
борьбы \cite{kar,pmixsin}.

Автор благодарит Леонида Посицельского за помощь в редактирова\-нии статьи.

\section{Постановка задачи}
Рассматривается игра двух лиц с нулевой суммой следующего вида.
Игроки обладают ресурсами $a\ge 0$ и $m\ge0$ ($m\in \Z$),
которые используют в промежутке времени $[0,1]$.
Первый игрок расходует свой ресурс непрерывным образом,
а второй --- единичными порциями.
Эффективность использования ресурса $j$-м игроком хaрактеризуется функцией
$P_j(t)$, $j=1,2$, равной вероятности достижения успеха при использовании
в момент $t$ единичного ресурса.
Функции $P_j(t)$ называются {\it функциями меткости}. Они непрерывно 
дифференцируемы, возрастают,
$P_j(0)=0$, $P_j(1)=1$, $P_j(t)<1$ при $t<1$.
Положим
$p(t)=1-P_1(t)$, $q(t)=1-P_2(t)$, $P(t)=(P_1(t), P_2(t))$.
Вектор-функция $P(t)$ называется {\it вектор-функцией эффективности}.
Вероятность $G(t,\D\g)$ достижения успеха при использовании
в момент $t$ ресурса $\D\g\ge0$ с функцией эффективности $P_j(t)$
вычисляется по формуле~\cite{psi,lk}:
\begin{equation}
G(t,\D\g)=1-(1-P_j(t))^{\D\g},\quad \D\g>0;\quad  G(t,0)=0.
\label{I2}
\end{equation}

Выигрыш $j$-го игрока в случае его успеха составляет $A_j$.
Пару $(A_1,A_2)$ обозначим $A$.
Выигрыш игроков равен $0$, если ни один из них не добился успеха
или они добились успеха одновременно.

 Пусть  $\a (t)$, $n(t)$ --- ресурсы игроков, оставшиеся к моменту $t$.
 Функции $\a(t)$, $n(t)$ назовем {\it функциями расхода ресурса}.
 Они не возрастают, функция $n(t)$ кусочно-постоянна и непрерывна слева.
 Пара функций $(\a(t), n(t))$ называется {\it партией}.
 Будем предполагать, что функция $\a(t)$ непрерывна на $[0,1]$ и
кусочно-непрерывно дифференцируема на $(0,1)$.
 Положим $\xi(t)=-\a'(t)$ и назовем функцию $\xi(t)$
{\it интенсивностью} расхода ресурса.

{\it Моментом действия} игрока называется момент времени, когда игрок
использует свой ресурс. У пулеметчика это точка убывания его функции
расхода ресурса ($\xi(t+0)>0$), у снайпера --- точка скачка этой функции.
Пусть  $\eta_i$, $i=1,2,\dots,m$
($0\le\eta_m\le\eta_{m-1}\le\dots\le\eta_1\le 1$) ---
моменты действия 2-го игрока. Вектор $\eta=(\eta_1,\,\eta_2,\dots\eta_m)$
назовем {\it вектором моментов действия} второго игрока.

Предположим, что на отрезке $[t_1,t_2]$ второй игрок не использует свой
ресурс. Вычислим вероятность того, что первый игрок
достигнет успеха на этом отрезке при интенсивности расхода $\xi(t)$.
В силу \eqref{I2} имеем:
\begin{align}
\nonumber
\varphi(\xi, t_1, t_2)&=
\lim_{N\to\infty}\lt(1-\prod_{i=1}^N\lt(1-P_1(\tau_i)\rt)^
{\xi(\tau_i)\D\tau_i}\rt)=\\
\label{phi}
&=1-\exp\int_{t_1}^{t_2}\xi(\tau)\ln (1-P_1(\tau))\,d\tau,
\end{align}
где $\xi(t)=-\alpha'(t)$.
Вероятность успеха первого игрока на отрезке $[t_1,t_2]$ выражается через
функцию расхода ресурса $\alpha(t)$ следующим образом
\begin{align}
\label{phi1}
\varphi(\alpha, t_1, t_2)=
1-\exp\left(-\int_{t_1}^{t_2}\ln (1-P_1(\tau))\,d\alpha(\tau)\right).
\end{align}
Обозначим $K(\alpha; \eta_1, \eta_2,\dots, \eta_m)$ математическое ожидание
выигрыша, получаемого 1-м игроком, если он расходует ресурс
в соответствии с функцией расхода ресурса $\alpha(t)$
в то время как 2-й  действует в моменты времени $\eta_k$ ($1\le k\le m$).
Тогда $K=0$ при $m=0$, $a=0$ и $K=A_1$ при $m=0$, $a>0$.
При $m\ge 1$  значение функции $K$ вычисляется по рекуррентной формуле:
\begin{align}
\label{defK}
K(\alpha; \eta_1, \eta_2,\dots, \eta_m)&=
A_1\varphi(\alpha, 0, \eta_m)-A_2(1-\varphi(\alpha, 0, \eta_m))P_2(\eta_m)+\\
&+(1-P_2(\eta_m))(1-\varphi(\alpha, 0, \eta_m))
K(\alpha_m; \eta_1, \eta_2,\dots, \eta_{m-1}), \nonumber
\end{align}
где
$$
\alpha_m(t)=
  \begin{cases}
    \alpha(\eta_m), &t\in[0,\eta_m)\\
    \alpha(t), &t\in[\eta_m,1].
  \end{cases}
$$

 Рассматриваемая игра, называемая шумной дуэлью,  представляет собой
модель конкурентной борьбы в условиях полной информации.
 В этой игре каждому игроку  известны ресурсы обоих соперников
в каждый момент времени,  и игроки непрерывно корректируют свое поведение
на основании поступающей информации о значении ресурса противника.
Cтратегия  первого игрока есть функция $\xi=u(t,\a,n)$,
определяющая интенсивность $\xi$ расхода ресурса в момент $t$
в зависимости от текущих значений $\a$ и $n$ ресурсов игроков.
 Стратегия 2-го игрока есть функция $\eta_n=v(\a,n)$, определяющая
очередной момент его действия по значениям  $\a$ и $n$ ресурсов.
 Мы будем задавать стратегии игроков только на том промежутке времени,
где $\alpha(t)n(t)>0$, предполагая что если один из игроков еще не
исчерпал свой ресурс, то он использует его так, чтобы обеспечить
единичную вероятность успеха.
 В силу условия $P_j(1)=1$  и формулы \eqref{phi1} это всегда возможно.
 Платежная функция этой игры есть функция $K(\xi,\eta)$, определенная
формулой \eqref{defK}, где $\xi$ и $\eta$ --- реализованные в ходе игры
интенсивность расхода ресурса 1-го игрока и вектор моментов действия 2-го.
 Описанную дуэль обозначим $G_{am}(P, A)$.

\section {$T$-партии и $T$-стратегии}
Обозначим через $\cal T$ множество  функциональных последовательностей
$$T(x)=(T_1(x), T_2(x), \dots, T_k(x), \dots),$$
удовлетворяющих следующим условиям:
\begin{enumerate}
\item[1.]
Функции $T_k(x)$ ($k\in\N$) определены и непрерывны на полупрямой
$[0,+\infty)$, непрерывно дифференцируемы на $(0,+\infty)$.
\item[2.]
$0 < T_k(x)\le 1$  ($x\ge 0$, $k\in\N$).
\item[3.]
$T'_k(x)<0$, $T_{k+1}(x)<T_{k}(x)$ ($x>0$, $k\in\N$).
\item[4.]
$T_k(0)=1$ ($k\in\N$).
\end{enumerate}
Пусть $T\in \cal T$. Партию $(\alpha, n)$ назовем {\it $T$-партией}, если
при $\alpha(t)n(t)>0$  выполняется неравенство
\begin{align}
t\le T_{n(t)}(\alpha(t)),
\label{Tan}
\end{align}
причем в тех точках, которые являются моментами действия,
неравенство выполняется как равенство.

Каждая последовательность $T\in \cal T$ определяет множество $T$-партий,
отличающиxся тем, кто из игроков использует ресурс в момент действия,
предписываемый последовательностью $T$
(возможны и одновременные действия игроков).
{\it Простейшими}  $T$-партиями назовем $T$-партии $(\alpha_j, n_j)$
($j=1,2$), в которых $j$-й игрок приступает к использованию своего ресурса
после того, как ресурс противника исчерпался. Функции расхода ресурса,
составляющие простейшие $T$-партии, имеют следующий вид:
\begin{align}
\label{baseT1}
&\alpha_1(t)=\begin{cases}
             a,\text{ при }t\in[0,T_1(a)];\\
             \text {обеспечивает расходимость }
            \Ds\int\limits_{T_1(a)}^1\ln(1-P_1(t))\,d\alpha_1(t),
            \end{cases}\\
\label{baseT2}
&n_1(t)=\begin{cases}
            m, &\text{при }t\in [0,T_m(a)];\\
            i, &\text{при }t\in (T_{i+1}(a),T_i(a)],\; 1\le i\le m-1;\\
            0, &\text{при }t\in (T_1(a),1],
        \end{cases}
\\
&\alpha_2(t)=\begin{cases}
              a,& \text{при }t\in[0,T_m(a)];\\
              T^{-1}_m(t),&\text{при }t\in(T_m(a),1],
             \end{cases}\\
&n_2(t)= m,\text{ при }t\in [0,1]\;(\eta_1=\eta_2=\dots=\eta_m=1).
\end{align}

\begin{lm}
\label{l:Tan1}
Пусть $T\in \cal T$. Значения платежной функции во всех  $T$-партиях дуэли
$G_{am}(P, A)$ совпадают в том и только том случае, когда выполнено
равенство
\begin{align}
\label{eqT1}
&\exp\lt(\int_0^x\ln(1-P_1(T_k(\a)))\,d\a\rt)+
\prod_{i=1}^k \lt(1-P_2(T_i(x))\rt)=1\\
&\text{ при любых }0<x\le a,\; 1\le k\le m.
\nonumber
\end{align}
При этом общее значение платежной функции во всех  $T$-партиях игры равно
\begin{align}
\label{vam1}
v_m(a)&=
A_1-(A_1+A_2)\exp\lt(\int_0^a\ln(1-P_1(T_m(\a)))\,d\a\rt)=\\
&=(A_1+A_2)\prod_{i=1}^m(1-P_2(T_i(a))-A_2.
\nonumber
\end{align}
\end{lm}

\begin{pf} \underline{Heoбходимость}.
Пусть $T\in \cal T$ и значения платежной функции во всех  $T$-партиях дуэли
$G_{am}(P, A)$ совпадают. Зафиксируем $x$ ($0< x \le a$), $k$ ($1\le k \le m$)
и положим $t^*=T_k(x)$.
Обозначим через $(\overline{\alpha}_j, \overline{n}_j)$ ($j=1,2$)
простейшие  $T$-партии игры $G_{xk}(P,A)$ и рассмотрим две $T$-партии
$(\alpha_j,n_j)$ ($j=1,2$) игры $G_{am}(P,A)$, удовлетворяющие следующим
условиям:
\begin{align}
&\alpha_1(t)=\alpha_2(t);\;n_1(t)=n_2(t)
               \text{ при }t\in[0,t^*];
\label{Tan1}
\\
&\alpha_j(t^*)=x;\; n_j(t^*)=k,\;j=1,2;
\label{Tan2}
\\
&\alpha_j(t)=\overline{\alpha}_j(t);\;
n_j(t)=\overline{n}_j(t)\text{ при }t\in(t^*,1],\;j=1,2.
\label{Tan3}
\end{align}

Обозначим $\overline{\eta}^j$ --- векторы моментов действия снайпера
в партиях $(\overline{\alpha}_j,\overline{n}_j)$ игры $G_{xk}(P,A)$
при $j=1,2$.
Вычислим значения платежной функции в этих партиях:
\begin{align}
\nonumber
K(\overline{\alpha}_1,\overline{\eta}_1)&=-A_2
\left(1-\prod_{i=1}^k \lt(1-P_2(T_i(x))\rt)\rt)
+A_1\prod_{i=1}^k \lt(1-P_2(T_i(x))\rt)=\\
\label{Tan4}
&=(A_1+A_2)\prod_{i=1}^k \lt(1-P_2(T_i(x))\rt)-A_2;\\
\nonumber
K(\overline{\alpha}_2,\overline{\eta}_2)&=A_1
\left(1-\exp\lt(\int_0^x\ln(1-P_1(T_k(\a)))\,d\a\rt)\rt)-\\
\nonumber
&-A_2\exp\lt(\int_0^x\ln(1-P_1(T_k(\a)))\,d\a\rt)=\\
\label{Tan5}
&=A_1-(A_1+A_2)\exp\lt(\int_0^x\ln(1-P_1(T_k(\a)))\,d\a\rt).
\end{align}
   По условию теоремы значения платежной функции игры $G_{am}(P,A)$ в партиях
$(\alpha_1,n_1)$ и $(\alpha_2,n_2)$ совпадают. Поэтому из условий
\eqref{Tan1}--\eqref{Tan3} следует равенство значений платежной функции
игры $G_{xk}(P,A)$ в партиях $(\overline{\alpha}_1,\overline{n}_1)$ и
$(\overline{\alpha}_2,\overline{n}_2)$.
Приравнивая
\eqref{Tan4} и \eqref{Tan5}, получим \eqref{eqT1}.

\underline{Достаточность}. Пусть $T\in{\cal T}$, равенство \eqref{eqT1}
 выполнено при всех $0<x\le a$,\; $1\le k\le m$ и $(\alpha, n)$ ---
произвольная $T$-партия. Покажем, что $K(\alpha, n)=v_m(a)$. Доказательство
проведем индукцией по числу единиц ресурса снайпера. При $m=0$ утверждение
леммы выполнено, так как $K=A_1$. Предположим, что равенство выполнено при
$n\le k-1$ и докажем его для $n=k$.
Пусть $\eta$ --- вектор  моментов действия снайпера в партии $(\alpha, n)$.
 Положим
$$\alpha_k=\alpha(\eta_k),\quad \alpha^k(t)=\min\{\alpha_k,\alpha(t)\}.$$
 Тогда
\begin{align}
\label{Tan6}
&K(\alpha, \eta_1,\dots,\eta_k)=A_1\varphi(\alpha, 0,\eta_k)-A_2P_2(\eta_k)
\left(1-\varphi(\alpha, 0,\eta_k)\rt)+\\
\nonumber
&+(1-P_2(\eta_k))\left(1-\varphi(\alpha, 0,\eta_k)\rt)
K(\alpha^k, \eta_1,\dots,\eta_{k-1}),
\end{align}
где
$\varphi(\alpha, 0,\eta_k)$
 --- вероятность того, что пулеметчик достиг успеха
на промежутке $[0,\eta_k)$. Согласно формуле  \eqref{phi1} имеем:
\begin{align}
\label{Tan7}
\varphi(\alpha, 0,\eta_k)=
1-\exp\lt(\int_{\alpha_k}^a\ln(1-P_1(T_k(\a)))\,d\a\rt).
\end{align}
По индукционному предположению
\begin{align}
\label{Tan8}
K(\alpha^k, \eta_1,\dots,\eta_{k-1})&=
(A_1+A_2)\prod_{i=1}^{k-1}(1-P_2(T_i(\alpha_k))-A_2.
\end{align}
Подставляя \eqref{Tan7} и \eqref{Tan8} в \eqref{Tan6}, получаем
\begin{align}
\label{Tan9}
&K(\alpha, \eta_1,\dots,\eta_k)=A_1+\\
\nonumber
&+(A_1+A_2)\exp\left(\int_{\alpha_k}^a\ln(1-P_1(T_k(\a)))\,d\a\right)
\left(
\prod_{i=1}^k(1-P_2(T_i(\alpha_k))-1\rt).
\end{align}
В силу \eqref{eqT1} имеем:
\begin{align}
\label{Tan10}
1-\prod_{i=1}^k\lt(1-P_2(T_i(\alpha_k))\rt)=
\exp\int_0^{\alpha_k}\ln(1-P_1(T_k(\a)))\,d\a.
\end{align}
Из \eqref{Tan9} с учетом \eqref{Tan10} окончательно получаем:
\begin{align*}
K(\alpha, \eta_1,\dots,\eta_k)=A_1-
(A_1+A_2)\exp\int_0^a\ln(1-P_1(T_k(\a)))\,d\a=v_k(a).
\end{align*}
Утверждение леммы по индукции доказано.
\end{pf}

Пусть $\{T_k(x)\}\in {\cal T}$.
{\it $T$-стратегиями\/} игроков назовем стратегии вида
\begin{align*}
&\xi^T(t)=
\begin{cases}
              0,& t<T_n(\a),\\
              -1/T'_n(\a),& t=T_n(\a);
\end{cases}
\quad\eta^T=T_n(\a),
\end{align*}
где $\a$, $n$ --- ресурсы игроков, оставшиеся к моменту $t$.

\begin{thm}
\label{t:Tan1}
Если последовательность $\{T_k(x)\}\in {\cal T}$ удовлетворяет
соотношению \eqref{eqT1} при всех $0<x\le a,$\; $1\le k\le m$,
то пара $T$-стратегий образует ситуацию равновесия игры $G_{am}(P,A)$.
Цена игры вычисляется по формуле
\eqref{vam1}.
\end{thm}

\begin{pf}
Пусть $\eta$ --- произвольный вектор моментов действия снайпера.
Предположим, что пулеметчик действует в соответствии с $T$-стратегией.
Его функция расхода ресурса, соответствующая $\eta$ при $T$-стратегии имеет
 вид:
\begin{align}
\label{Tan1.1}
&\alpha^T(t)=
\begin{cases}
              \alpha_{k+1},& t\in (\eta_{k+1},T_k(\a_{k+1}));\\
              \alpha^T_k(t),& t\in [T_k(\a_{k+1}),\eta_k],
                    \text{ если } \eta_k>T_k(\a_{k+1}),
\end{cases}
\end{align}
где  $\alpha^T_k(t)$ --- функция, определенная на отрезке
$[T_k(\alpha_{k+1}),1]$ и обратная к $T_k(x)$, продолженная на отрезок
$[0,T_k(\alpha_{k+1})]$ как константа $\alpha_{k+1}$,
$\alpha_k=\alpha^T_k(\eta_k)$ ($k=1,2\dots,m$),
$\eta_{m+1}=0$, $\alpha_{m+1}=a$.

Покажем, что если $T_k(x)$ удовлетворяет  \eqref{eqT1} при всех
$0<x\le a,$\; $1\le k\le m$, то выполняется неравенство
$$K(\alpha^T;\eta)\ge v_m(a).
$$

Прежде всего отметим, что при использовании пулеметчиком $T$-стратегии
неравенство $\eta_k>T_k(\alpha_k)$ невозможно, так как, начиная с момента
$T_k(\alpha_{k+1})$, пулеметчик начинает расходовать ресурс в соответствии
с функцией $\alpha^T_k(t)$, обеспечивающей равенство $t=T_k(\alpha^T(t))$,
и, в результате этого, очередной момент действия снайпера,
предписываемый $T$-статегией,
отодвигается.

Если при всех  $1\le k\le m$ $\eta_k=T_k(\alpha)$, то реализуется простейшая
$T$-партия \eqref{baseT1} и, согласно лемме~\ref{l:Tan1}, выполняется
равенство $K(\alpha^T;\eta)= v_m(a)$.

Если снайпер расходует ресурс, не дожидаясь наступления момента действия
$T_k(\alpha^T(t))$, то найдутся два таких номера $k$ и $l$
($1\le l\le k\le m$), что
\begin{align}
\label{Tan1.1a}
&\eta_i=T_i(\alpha_{i}) \text{ при } k+1\le i\le m;\quad
\eta_{l-1}=T_{l-1}(\alpha_{l-1});\\
\label{Tan1.1b}
&\eta_i<T_i(\alpha_{i}) \text{ при } l\le i\le k.
\end{align}
Тогда по определению $T$-стратегии ресурс пулеметчика на промежутке
$(\eta_{k+1},T_l(\a_{k+1}))$ не расходуется: $\alpha^T(t)=\alpha_{k+1}$.
Определим вектор  $\eta^1$ следующим образом:
\begin{align*}
&\eta^1_i=
\begin{cases}
                   \eta_i,& i=1,2,\dots,l-1,k+1,k+2,\dots,m;\\
                   T_i(\a_{k+1}),& i=l,l+1,\dots,k.
\end{cases}
\end{align*}
Вычислим $K(\alpha^T;\eta)$ и $K(\alpha^T;\eta^1)$, представив платежную
функцию в виде суммы трех слагаемых, соответствующих промежуткам $[0,\eta_k)$,
 $[\eta_k, T_l(\a_{k+1}))$, $[T_l(\a_{k+1}),1]$.
Получим:
\begin{align}
\label{Tan1.2}
&K(\alpha^T,\eta)=K_{[0,\eta_k)}-A_2\Psi(\eta_k)+
\Psi(\eta_k)\prod_{i=l}^kq(\eta_i)(A_2+K_{[T_l,1]});\\
\label{Tan1.3}
&K(\alpha^T,\eta^1)=K_{[0,\eta_k)}-A_2\Psi(\eta_k)+
\Psi(\eta_k)\prod_{i=l}^kq(T_i(\alpha_{k+1}))(A_2+K_{[T_l,1]}),
\end{align}
где  $K_{[0,\eta_k)}$ и $K_{[T_l,1]}$ --- математические ожидания выигрыша
первого игрока на промежутках $[0,\eta_k)$ и $[T_l(\a_{k+1}),1]$ при
функции расхода ресурса пулеметчика $\alpha^T$ и векторе $\eta$ моментов
действия снайпера, а $\Psi(\eta_k)$ --- вероятность того, что при указанных
функции расхода ресурса пулеметчика и векторе моментов действия снайпера
игроки не достигли успеха к моменту $\eta_k$.
Поскольку функция $q(t)$ убывает, то сравнивая
\eqref{Tan1.2} с \eqref{Tan1.3} и учитывая \eqref{Tan1.1a}, \eqref{Tan1.1b},
получаем  неравенство
$$K(\alpha^T,\eta)\ge K(\alpha^T,\eta^1).$$
Повторяя  описанную процедуру, построим $r$ векторов
$\eta^1,\eta^2,\dots,\eta^r$ ($r<m$) таких, что
\begin{align*}
&K(\alpha^T,\eta)\ge K(\alpha^T,\eta^1)\ge K(\alpha^T,\eta^2)
\ge\dots\ge K(\alpha^T,\eta^r) \\
&\text {и }
\eta^r_i=T_i(\alpha_i) \text{ при всех } 1\le i\le m.
\end{align*}
По лемме~\ref{l:Tan1} $K(\alpha^T;\eta^r)= v_m(a)$ и, следовательно,
$K(\alpha^T;\eta)\ge v_m(a)$.

Пусть теперь  $\alpha(t)$ --- произвольная функция расхода ресурса пулеметчика,
 а $\eta^T$ --- реализация вектора моментов действия снайпера, соответствующая
$\alpha(t)$ при $T$-стратегии. Покажем, что
$K(\alpha;\eta^T)\le v_m(a)$.  Обозначим через $\alpha^T$ реализацию функции
расхода ресурса пулеметчика, соответствующую $\eta^T$ при $T$-стратегии
($\alpha^T$ находится по формуле \eqref{Tan1.1}). Если при всех $t\in[0,1]$
таких, что $\alpha(t)n(t)>0$, функция $\alpha(t)$ совпадает с $\alpha^T(t)$,
то реализуется $T$-партия, и в силу леммы~\ref{l:Tan1}
$K(\alpha;\eta^T)= v_m(a)$. В противном случае найдется $t_*\in[0,1]$ такое,
что $\alpha(t_*)=\alpha^T(t_*)$
и для некоторого $\eps>0$ при всех $t\in(t_*,t_*+\eps)$
выполнено неравенство $\alpha(t)<\alpha^T(t)$. Противоположное неравенство
невозможно, поскольку снайпер использует $T$-стратегию.
Положим
$$
t^*=\sup \{t: \alpha(t')<\alpha^T(t') \text{ при всех } t'\in(t_*,t)\}.
$$
На отрезке $[0,1]$ содержится не более, чем счетное множество отрезков вида
$[t_*,t^*]$. Перенумеруем их, и пусть $[t_1, t'_1]$ --- первый тaкой отрезок.
Определим  $\alpha_1(t)$ следующим образом:
\begin{align*}
&\alpha_1(t)=
\begin{cases}
             \alpha^T(t), &\text {при } t\in(t_1, t'_1);\\
             \alpha(t),   &\text {при } t\notin(t_1, t'_1).
\end{cases}
\end{align*}
Пусть $t_1\in[\eta^T_k,\eta^T_{k-1}]$, тогда по определению $T$-стратегии
 снайпера $t'_1\in[\eta^T_k,\eta^T_{k-1}]$.
Вычислим $K(\alpha;\eta^T)$ и $K(\alpha_1;\eta^T)$, представив платежную
функцию в виде суммы трех слагаемых, соответствующих промежуткам $[0,t_1)$,
 $[t_1, t'_1)$, $[t'_1,1]$.
Получим:
\begin{align}
\label{Tan1.4}
K(\alpha,\eta^T)&=K_{[0,t_1)}+A_1\Psi(t_1)
-\Psi(t_1)(A_1-
K_{[t_1',1]})\exp\int\limits_{t_1}^{t'_1}\mu(t)\,d\alpha(t);\\
\label{Tan1.5}
K(\alpha_1,\eta^T)&=K_{[0,t_1)}+A_1\Psi(t_1)
-\Psi(t_1)(A_1-
K_{[t_1',1]})\exp\int\limits_{t_1}^{t'_1}\mu(t)\,d\alpha^T(t),
\end{align}
где  $K_{[0,t_1)}$ и $K_{[t_1',1]}$ --- математические ожидания выигрыша
первого игрока на промежутках $[0,t_1)$ и $[t_1',1]$ при
функции расхода
ресурса пулеметчика $\alpha(t)$ и векторе $\eta^T$ моментов действия снайпера,
 $\Psi(t_1)$ --- вероятность того, что при указанных
функции расхода
ресурса пулеметчика и векторе моментов действия снайпера игроки не достигли
успеха к моменту $t_1$, $\mu(t)=-\ln p(t)$.
Заметим, что
\begin{align}
\label{Tan1.6}
\int\limits_{t_1}^{t'_1}\mu(t)\,d\alpha(t)>
\int\limits_{t_1}^{t'_1}\mu(t)\,d\alpha^T(t)
\end{align}
Действительно, интегрируя по частям, получим
\begin{align}
\label{Tan1.7}
&\int_{t_1}^{t'_1}\mu(t)\,d\alpha(t)
=\mu(t'_1)\int_{t_1}^{t'_1}\,d\alpha(t)-
\int_{t_1}^{t'_1}
\left(\int_{t_1}^{t}\,d\alpha(\tau)\right)\,d\mu_1(t),\\
\label{Tan1.8}
&\int_{t_1}^{t'_1}\mu(t)\,d\alpha^T(t)
=\mu(t'_1)\int_{t_1}^{t'_1}\,d\alpha^T(t)-
\int_{t_1}^{t'_1}
\left(\int_{t_1}^{t}\,d\alpha^T(\tau)\right)\,d\mu_1(t).
\end{align}
По построению отрезка $[t_1,t'_1]$, при всех $t\in[t_1,t'_1]$ выполняется
неравенство
\begin{align}
\label{Tan1.9}
           \int_{t_1}^{t}\,d\alpha(\tau)\le
           \int_{t_1}^{t}\,d\alpha^T(\tau),
\end{align}
причем равенство достигается только при $t=t'_1$. Из
\eqref{Tan1.7} и \eqref{Tan1.8} в силу \eqref{Tan1.9} вытекает
\eqref{Tan1.6}. Сравнивая $K(\alpha; \eta^T)$ с $K(\alpha_1; \eta^T)$
(формулы \eqref{Tan1.4},  \eqref{Tan1.5}) и
учитывая \eqref{Tan1.6}, получаем, что
\begin{align*}
 K(\alpha,\eta^T)\le K(\alpha_1,\eta^T).
\end{align*}

Повторяя описанную процедуру, строим последовательность функ\-ций $\alpha_k$
такую, что
\begin{align*}
 K(\alpha_k,\eta^T)\le K(\alpha_{k+1},\eta^T),\;k\in\N.
\end{align*}
Таким образом,
\begin{align}
\label{Tan1.10}
 K(\alpha,\eta^T)\le K(\alpha_{k},\eta^T), \text{ при каждом }k\in\N.
\end{align}
Положим $\alpha^*(t)=\Ds\lim_{k\to\infty}\alpha_k(t)$. Переходя к пределу при
$k\to\infty$ в неравенстве
\eqref{Tan1.10} в силу 1-й теоремы Хелли \cite{kf} заключаем, что
\begin{align*}
 K(\alpha,\eta^T)\le K(\alpha^*,\eta^T).
\end{align*}
Поскольку $\alpha^*(t)=\alpha^T(t)$ при всех  $t\in\{t:n(t)>0\}$, то согласно
лемме~\ref{l:Tan1} $K(\alpha^*;\eta^T)= v_m(a)$ и, следовательно,
$K(\alpha;\eta^T)\le v_m(a)$.
\end{pf}

\begin{cor}
\label{c:Tan1}
Если функция $P_2(t)$ строго возрастает на отрезке $[0,1]$, то существует
не более одной последовательности $\{T_k(x)\}\in {\cal T}$,
удовлетворяющей \eqref{eqT1} при всех $x\ge 0$, $k\in\N$.
\end{cor}

\begin{pf}
Предположим, существуют две последовательности $T^1,\;T^2\in{\cal T}$,
удовлетворяющие \eqref{eqT1} при всех $x\ge 0$, $k\in\N$.
Пусть
\begin{align}
\label{cor.1}
l=\min\{k:T_k^1\neq T_k^2\} \quad\text{и}\quad
T_l^1(a)\neq T_l^2(a),\; a>0.
\end{align}
В силу теоремы~\ref{t:Tan1} игра $G_{al}(P,A)$ имеет цену, равную
значению платежной функции в  $T^1$- и $T^2$-партиях, т.~е.\ выполнено
равенство
\begin{align}
\label{cor.2}
(A_1+A_2)\prod_{i=1}^{l}\lt(1-P_2(T^1_i(a))\rt)-A_2=
(A_1+A_2)\prod_{i=1}^{l}\lt(1-P_2(T^2_i(a))\rt)-A_2.
\end{align}
Но по определению $l$ при всех $i<l$ выполнено равенство $T_i^1(a)=T_i^2(a)$,
поэтому из \eqref{cor.2} с учетом строгой монотонности функции $P_2(t)$
заключаем, что $T_l^1(a)=T_l^2(a)$.
Полученное противоречие доказывает единственность последовательности
$\{T_k(x)\}$.
\end{pf}

Пусть последовательность $\{T_k(x)\}\in {\cal T}$ удовлетворяет
соотношению \eqref{eqT1} при всех $a\ge 0$ и $k\in \N$.
Введем обозначения:
\begin{align}
\label{eqT2.1a}
\pi_0(x)=1;\;
 \pi_{k}(x)=q(T_{k}(x))\pi_{k-1}(x)\;(k\in\N).
\end{align}
Тогда
\begin{align}
\exp\lt(\int_0^x\ln(p(T_k(\a))\,d\a\rt)=1-\pi_{k}(x)
\label{eqT2.1}
\end{align}
Дифференцируя \eqref{eqT2.1} по  $x$, получим
\begin{align}
\label{eqT2.2}
\ln(p(T_k(x))=-\frac{\pi'_{k}(x)}{1-\pi_{k}(x)}
\end{align}
Напишем рекуррентную формулу для $\pi'_{k}(x)$:
\begin{align}
\label{eqT2.3}
\pi'_{k}(x)=q'(T_k(x))T'_k(x)\pi_{k-1}(x)+q(T_k(x))\pi'_{k-1}(x).
\end{align}
Из \eqref{eqT2.2} и \eqref{eqT2.3} следует, что последовательность
$\{T_k(x)\}$ удовлетворяет системе обыкновенных дифференциальных уравнений
\begin{align}
\label{difeq}
&\frac{dT_k}{dx}=\phi(T_1,T_2, \dots, T_k), \quad
x\ge 0;\;  k\in \N,
\end{align}
где
\begin{align*}
&\phi(T_1,T_2, \dots, T_k)=\\
&=-\frac{\left(
1-\Ds\prod_{i=1}^k q(T_i)\rt)\ln p(T_k)-q(T_k)
\lt(1-\Ds\prod_{i=1}^{k-1}q(T_i)\rt)\ln p(T_{k-1})}
{q'(T_k)\Ds\prod_{i=1}^{k-1}q(T_i)}.
\end{align*}

\begin{lm}
\label{l:Tan2}
Пусть $p(t)$, $q(t)$ непрерывно дифференцируемы на $(0,1]$, $p(0)=q(0)=1$;
$p(1)=q(1)=0$; $p(t)>0$ при $t<1$; $p'(t)\le 0$; $q'(t)<0$.
Тогда система обыкновенных дифференциальных уравнений  \eqref{difeq}
с начальными условиями
\begin{align}
T_k(0)=1 \quad (k\in \N)
\label{difeq1}
\end{align}
имеет решение на полупрямой $x>0$, причем $\{T_k(x)\}\in {\cal T}$.
\end{lm}

\begin{pf}
Доказательство проведем по индукции. Первое уравнение системы \eqref{difeq}
имеет вид:
\begin{align}
\label{difeq2}
&\frac{dT_1}{dx}=
-\frac{(1-q(T_1))\ln p(T_1)}{q'(T_1)}.
\end{align}
Интегрируем \eqref{difeq2} с начальным условием $T_1(0)=1$. Получаем
\begin{align}
\label{difeq3}
x(T_1)=
\int\limits_{T_1}^1\frac{q'(\tau)\,d\tau}{(1-q(\tau))\ln p(\tau)}.
\end{align}
Функция $T_1(x)$ является обратной к $x(T_1)$. Проверим, что она
удовлетворяет условиям 1--4.

Покажем, что $T_1(x)$ определена на полупрямой $[0,+\infty)$.
Выберем $\delta>0$ такое, что $-\ln p(\delta)<1$. Тогда
$$
x(t)\ge
\int\limits_{t}^{\delta}\frac{(1-q(\tau))'\,d\tau}{(1-q(\tau))}+x(\delta)=
\ln(1-q(\delta))-\ln(1-q(t))+x(\delta).
$$
Поскольку правая часть неравенства стремится к $+\infty$ при $t\to +0$, то
при $x\to +\infty$  $T_1\to +0$ и, следовательно, функция $T_1(x)$ определена
при всех $x>0$ и, согласно начальному условию,
$T_1(0)=1$.

В силу \eqref{difeq2} $T'_1(x)<0$. Таким образом, $T_1(x)$ убывает от 1 при
$x=0$ к $0$ при $x\to +\infty$. Следовательно,
$0 < T_1(x)\le 1$.

Предположим, что при $1\le i\le k-1$ решения $T_i(x)$ системы \eqref{difeq},
\eqref{difeq1} существуют и удовлетворяет условиям 1--4. Подставив их в $k$-е
уравнение системы \eqref{difeq}, получим
\begin{align}
\label{difeq4}
\frac{dT_k}{dx}=\Phi_k(T_k,x),
\end{align}
где
\begin{align}
\label{difeq5}
\Phi_k(t,x)=\phi_k(T_1(x),T_2(x),\dots, T_{k-1}(x),t).
\end{align}
Покажем, что уравнение
\eqref{difeq4} с начальным условием
\begin{align}
\label{difeq6}
T_k(0)=1
\end{align}
имеет решение $T_k(x)$ и для него выполняются условия монотонности
\begin{align}
\label{difeq7}
T'_k(x)<0,\; T_{k}(x)<T_{k-1}(x) \text{ при всех }x>0.
\end{align}
Заметим, что при всех $x>0$ справедливо неравенство:
\begin{align}
\label{difeq8}
\Phi_k(T_{k-1}(x),x)\le\frac{T'_{k-1}(x)}{q(T_{k-1}(x))}<T'_{k-1}(x)<0.
\end{align}
Действительно,
\begin{align}
\label{difeq9}
\Phi_k(T_{k-1}(x),x)=
-\frac{\left(1-q(T_{k-1}(x))\rt)\ln p(T_{k-1}(x))}
{q'(T_{k-1}(x))\pi_{k-1}(x)},
\end{align}
где $\pi_{k}(x)$ --- функции, определенные формулами \eqref{eqT2.1a}.
С другой стороны, по предположению индукции
$$T_{k-1}(x)<T_{k-2}(x),$$
поэтому из \eqref{difeq}
вытекает неравенство
\begin{align}
\label{difeq10}
T'_{k-1}(x)\ge
-\frac{\left(1-q(T_{k-1}(x))\rt)\ln p(T_{k-1}(x))}
{q'(T_{k-1}(x))\pi_{k-2}(x)}.
\end{align}
Сравнивая \eqref{difeq9}  и \eqref{difeq10}, получаем \eqref{difeq8}.
Обозначим числитель дроби, стоящей в правой части уравнения \eqref{difeq}
$F_{k}(t,x)$, т.~е.
\begin{align*}
F_{k}(t,x)=
(1-\pi_{k-1}(x)q(t))\ln p(t)-
(1-\pi_{k-1}(x))q(t)\ln p(T_{k-1}(x)).
\end{align*}
Покажем, что уравнение
\begin{align}
\label{difeq11}
F_{k}(t,x)=0,\;x>0,
\end{align}
задает неявную функцию $t=f_k(x)$, обладающую следующими свой\-ствами:
\begin{enumerate}
\item[1)]
функция $f_k(x)$ ($k\in\N$) определена и непрерывно дифференци\-руема
на полупрямой $(0,+\infty)$;
\item[2)]
$f_k(x)\to 1$
 при $x\to +0$;
\item[3)]
$f'_k(x)<0$\;  ($x>0$);
\item[4)]
 $f_k(x)<T_{k-1}(x)$ \;
 ($x>0$, $k\ge 2$).
\end{enumerate}
   Из \eqref{difeq8} вытекает, что
$$F_k(T_{k-1}(x),x)<0 \text{ при всех } x>0.$$
С другой стороны, для любого $x>0$
$$
\lim_{t\to +0}F_{k}(t,x)=
-(1-\pi_{k-1}(x))\ln p(T_{k-1}(x))>0.
$$
Таким образом, при любом $x>0$ уравнение \eqref{difeq11} имеет корень
$$t=f_k(x)\in(0,T_{k-1}(x)).$$
Для доказательства единственности корня проверим, что
$$\frac{\partial F_{k}}{\partial t}<0 \text { при всех } x>0,\;t\in(0,1).$$
Действительно,
\begin{align*}
\frac{\partial F_{k}}{\partial t}&=
\frac{p'(t)}{p(t)}
(1-\pi_{k-1}(x)q(t))-
q'(t)\pi_{k-1}(x)\ln p(t)
-\\&
-q'(t)\pi_{k-1}(x)\ln p(T_{k-1}(x))<0.
\end{align*}
Из теоремы о неявной функции следует, что уравнение
\eqref{difeq11} определяет неявную функцию $t=f_{k}(x)$, дифференцируемую
 на полуоси $x>0$. Покажем, что $f'_{k}(x)<0$. Поскольку по доказанному
$F'_{k}(t)<0$, то, в силу соотношения
$$ (f_k)'_x=-\frac{(F_{k})'_t}{(F_{k})'_x},
$$
достаточно установить, что
$$ (F_{k})'_x<0 \text{ при всех } x>0,\; t<T_{k-1}(x).
$$
Учитывая, что в силу \eqref{eqT2.2}
$$\pi'_{k-1}(x)=-(1-\pi_{k-1}(x))\ln(p(T_{k-1}(x)),
$$
имеем
\begin{align*}
 (F_{k})'_x&=
-(1-\pi_{k-1}(x))\ln(p(T_{k-1}(x))\left(\ln p(t)-\ln(p(T_{k-1}(x))\rt)-\\
&-\frac{p'(T_{k-1}(x))}{p(T_{k-1}(x))}T'_{k-1}(x)(1-\pi_{k-1}(x))q(t).
\end{align*}
Поскольку $t<T_{k-1}(x)$ и $T'_{k-1}(x)<0$, то
 $(F_{k})'_x<0$.

Покажем теперь, что $f_k(x)\to 1$ при $x\to +0$. По доказанному $f_k(x)$
монотонно убывает на полуоси $x>0$. Учитывая, что
$$f_k(x)<T_{k-1}(x)<1,$$
заключаем, что существует предел  $f_k(x)$ при $x\to +0$ и
$$\lim_{x\to +0}f_k(x)=c\le 1.$$
Предположим, что $c<1$. Подставляя $t=f_k(x)$ в \eqref{difeq11}, получаем:
\begin{align}
\label{difeq12}
(1-\pi_{k-1}(x)q(f_k(x)))\ln p(f_k(x))=
(1-\pi_{k-1}(x))q(f_k(x))\ln p(T_{k-1}(x)),
\end{align}
При $x\to +0$ правая часть равенства
\eqref{difeq12} стремится к $-\infty$, а предел левой части равен
$\ln p(c)>-\infty$. Из полученного противоречия следует, что
$$\lim_{x\to +0}f_k(x)=1.$$

Перейдем к построению функции $T_k(x)$. Сопоставим каждому $a>0$ решения
$y_a(x)$ и $z_a(x)$
 уравнения \eqref{difeq4} на полуоси $x\ge a$,
 удовлетворяющие начальным условиям
$$y_a(a)=T_{k-1}(a);\quad z_a(a)=f_{k}(a).
$$
Для этого уравнения в области $x>0$, $0<t<1$  выполнены условия теоремы
существования и единственности решения. В силу неравенств
\eqref{difeq8} и $f'_{k-1}(x)<0$, кривые
$y_a(x)$ и $z_a(x)$ при $x>a$ лежат между кривыми
$$t=f_{k}(x) \text{ и } t=T_{k-1}(x).
$$
Возьмем $c>0$ и положим $\Ds\inf_{a>0}y_a(c)=b$.
Пусть $\tilde{y}(x)$ --- решение уравнения   \eqref{difeq4}, удовлетворяющее
начальному условию $\tilde{y}(c)=b$. В силу теоремы единственности
интегральная кривая $\tilde{y}(x)$  при любом $a>0$ лежит строго между кривыми
$y_a(x)$ и $z_a(x)$ при всех $x>0$. Отсюда следует, что $\tilde{y}(x)$
можно продолжить на полуось $x\ge 0$, причем $\tilde{y}(0)=1$ и
$$f_{k}(x)<\tilde{y}(x)<T_{k-1}(x) \text{ при всех } x>0.
$$
Поскольку $F_k(t,x)$ убывает по $t$ при каждом $x>0$, то
$$F_k(\tilde{y}(x),x)<F_k(f_{k}(x),x)=0$$
и, следовательно, $\tilde{y}'(x)<0$.  Таким образом,
$\tilde{y}(x)$ удовлетворяет всем условиям, наложенным на функцию $T_k(x)$
и, следовательно, существование этой функции доказано.
Полагаем $T_k(x)=\tilde{y}(x)$.
\end{pf}

\begin{lm}
\label{l:Tan2a}
Пусть $p(t)$, $q(t)$ непрерывно дифференцируемы на $(0,1]$, $p(0)=q(0)=1$;
$p(1)=q(1)=1$; $p(t)>0$ при $t<1$; $p'(t)\le 0$; $q'(t)<0$.
Тогда решение системы дифференциальных уравнений  \eqref{difeq},
\eqref{difeq1} является единственным.
\end{lm}

\begin{pf}
Доказательство проведем по индукции. При $k=1$ функция $T_1(x)$ есть
обратная к $x(T_1)$, однозначно определяемой формулой
\eqref{difeq3}.
Предположим, что при $1\le i\le k-1$  система \eqref{difeq},
\eqref{difeq1} имеет единственное решение $T_i(x)$.
Докажем, что задача \eqref{difeq4},\eqref{difeq6}
имеет единственное решение. Пусть
$T_k(x)$ --- решение этой задачи,
построенное в лемме~\ref{l:Tan2}, а
$y(x)$ --- произвольное решение этой задачи.
Докажем, что
$$y(x)\equiv T_k(x).$$
Рассмотрим 2 случая.
\begin{enumerate}
\item[1.]
Интегральная кривая $y(x)$ уравнения \eqref{difeq4} при всех $x>0$ лежит
строго между графиками функций
$$y=f_{k}(x) \text{ и } y=T_{k-1}(x).
$$
Тогда, по доказанному, $y(x)$ удовлетворяет условиям
\eqref{difeq6}, \eqref{difeq7} и, в силу следствия~\ref{c:Tan1},
$$y(x)\equiv T_k(x).
$$
\item[2.]
При некотором $x>0$ одно из неравенств
$$f_{k}(x)<y(x)<T_{k-1}(x)
$$
нарушается. Покажем, что в этом случае кривая $y(x)$ не проходит через точку
$(0,1)$, т.~е.\ в этом случае начальное условие \eqref{difeq6} не выполняется.
\begin{enumerate}
\item[а)]
Пусть существует
$x_0>0$ такое, что
$y(x_0)<f_{k}(x_0)$.
Тогда
$y'(x)>0$ при всех $x\in(0,x_0)$ и, следовательно,
$y(x)<f_{k}(x_0)$ при всех $0<x<x_0$,
т.~е.\ $y(0)<1$.
\item[б)]
Пусть существует
$x_1>0$ такое, что
$y(x_1)>T_{k-1}(x_1)$.
Тогда  согласно неравенству
\eqref{difeq8}
$y(x)>T_{k-1}(x)$ при всех $0<x<x_1$.
Поскольку $q'(t)$ непрерывна и  $q'(t)<0$  при $t\in(0,1]$, то существуют
числа $c_1$ и $c_2$ такие, что
$$c_1\le -q'(t)\le c_2 \text { при } t\in [T_{k-1}(x_1),1].
$$
Из монотонного убывания функции $F_k(t,x)$ по $t$ и неравенства
\eqref{difeq8} следует, что при $t>T_{k-1}(x)$
\begin{align}
\label{difeq13}
\Phi_k(t,x)=-\frac{F_k(t,x)}{q'(t)\pi_{k-1}(x)}<
-\frac{F_k(T_{k-1}(x),x)}{q'(t)\pi_{k-1}(x)}<
\frac{T'_{k-1}(x)q'(T_{k-1}(x))}{q'(t)q(T_{k-1}(x))}.
\end{align}
Подберем такое $\delta>0$, что $q(T_{k-1}(\delta))<c_1/c_2$ и положим
$\delta_1=\min(x_1,\delta)$.
Тогда при $x\in(0,\delta)$, $t>T_{k-1}(x)$ имеет место неравенство
$$\Phi_k(t,x)<T'_{k-1}(x),$$
поэтому существует $\varepsilon>0$ такое, что при $x\in(0,\delta_1)$
выполняется неравенство $y(x)-T_{k-1}(x)>\varepsilon$, откуда следует, что
$y(x)$ не проходит через точку $(0,1)$.
\end{enumerate}
\end{enumerate}
\end{pf}

\begin{rem}
\label{Rem1}
Пусть $P_2(t)=t$, а $y(x)$ и $z(x)$ --- решения уравнения \eqref{difeq4}
на полуоси $[a,+\infty)$ с начальными условиями $y(a)=y_0$, $z(a)=z_0$, причем
$$f_k(a)\le z_0<y_0\le T_{k-1}(a).$$
Тогда разность $e(x)=y(x)-z(x)$ убывает по $x$.
\end{rem}

\begin{pf}
Рассмотрим производную разности:
$$
e'(x)=\Phi_k(y(x),x)-\Phi_k(z(x),x).
$$
Из соотношения
$$
\left(\Phi_k(t,x)\right)'_t=\left(F_k(t,x)\right)'_t/\pi_{k-1}(x),
$$
учитывая неравенство $(F_k)'_t<0$, установленное при доказательстве
леммы~\ref{l:Tan2}, получаем, что  $\left(\Phi_k\right)'_t<0$ при всех $x>0$.
Следовательно, $\Phi_k(t,x)$ убывает по $t$ при каждом $x>0$. Поскольку
$y(x)>z(x)$, то $e'(x)<0$ и, следовательно, разность $y(x)-z(x)$ убывает по $x$.
\end{pf}

Из теоремы~\ref{t:Tan1} и леммы~\ref{l:Tan2} следует
\begin{thm}
\label{t:Tan2}
Если в шумной дуэли пулеметчика со снайпером функции меткости игроков
$P_j(t)$ ($j=1,2$) непрерывно дифференцируемы на $(0,1]$,
$P_j(0)=0$, $P_j(1)=1$, $P_j(t)<1$ при $t<1$,  $P'_1(t)\ge 0$, $P'_2(t)> 0$,
то игра имеет ситуацию равновесия в чистых стратегиях.
Цена и оптимальные стратегии игры имеют вид:
\begin{align}
\label{fin}
&v_m(a)=(A_1+A_2)\prod_{i=1}^m(1-P_2(T_i(a))-A_2;\\
&\xi^T(t)=
\begin{cases}
              0,& t<T_n(\a),\\
              -1/T'_n(\a),& t=T_n(\a);
\end{cases}
\quad\eta^T=T_n(\a),
\end{align}
где $\a$, $n$ --- ресурсы игроков, оставшиеся к моменту $t$, а
$\{T_k(x)\}$ --- решение системы уравнений \eqref{difeq}, \eqref{difeq1}.
\end{thm}

\section{Выводы}

\begin{enumerate}
\item
    Оптимальные стратегии $\xi^T$, $\eta^T$ предписывают игрокам не
использовать свои ресурсы до момента времени $T_m(a)$ (назовем его моментом
начала расхода).
\item
Момент начала расхода есть функция  $T(\alpha,\beta)$ от текущих
значений  $\alpha$ и $\beta$ ресурсов игроков, не зависит от хода
игры до текущего момента и  является общим для обоих игроков.
\item
Чем больше ресурсы игроков, тем
раньше наступает момент начала расхода, так как $T_k(x)$ убывает по $x$ и $k$.
\item
При оптимальном поведении в момент начала расхода один из игроков (или оба)
начинают действовать.
\item
Если в момент начала расхода 1-й игрок начинает действовать в соответствии
с оптимальной функцией расхода, то он продолжает действовать вплоть до
очередного момента действия снайпера и в течение всего этого времени
выполняется равенство
\begin{equation}
t=T_{m(t)}(\a(t)).
\label{Concl}
\end{equation}
\item
Если в момент начала расхода 2-й игрок действует, то
его ресурс уменьшается на единицу и
игроки прекращают действовать до момента
начала расхода, соответствующего текущим значениям ресурсов.
\item
Если пулеметчик придерживается оптимальной стратегии, то оптимальная стратегия
снайпера предписывает  ему действовать в любой момент, когда
равенство \eqref{Concl} выполняется (возможно, одновременно с пулеметчиком,
прерывая его действия), или отложить свои действия до момента окончания игры
$t=1$, что не влияет на значение платежной функции.
\item
Если один из игроков использует $T$-стратегию, то реализуется партия,
удовлетворяющая условию $t\le T_{m(t)}(\a(t))$ при всех $t\in[0,1]$, 
таких что $\a(t)m(t)>0$.
\item
Если оба игрока действуют согласно $T$-стратегиям, то они расходуют ресурс
только в те моменты времени, когда выполнено равенство \eqref{Concl}.
При этом реализуется одна из $T$-партий. Во всех таких партиях
платежная функция принимает постоянное значение.
\end{enumerate}

Отметим, что последовательность функций $T_k(x)$ есть непрерывный аналог
бесконечной  матрицы $\{t_{mn}\}$ ($m\in\N$, $n\in\N$)
<<моментов наилучшего действия>> шумной дуэли с дискретным расходом
ресурса у обоих игроков \cite{fk}.

\section{Приложение. Численное решение игры}

Мы будем рассматривать случай $P_2(t)=t$, что не ограничивает общности.
Действительно, пусть $P_2(t)\not\equiv  t$, $P'_2(t)>0$ и $P_2(t)>0$ при
$t\in (0,1]$.
Тогда сделаем замену переменных $\tau=P_2(t)$ и будем решать игру
$G_{am}(\tilde P,A)$, где  $\tilde P=(P_1(P_2^{-1}(\tau)),\tau)$.
Очевидно, цены игр $G_{am}(P,A)$ и $G_{am}(\tilde P,A)$ совпадают,
а оптимальные $T$-стратегии игры $G_{am}(P,A)$ определяются последовательностью
$T_k(x)=P_2^{-1}(\tilde T_k(x))$, где $\tilde T_k(x)$ --- решение задачи
\eqref{difeq}, \eqref{difeq1} для игры $G_{am}(\tilde P,A)$.

В силу \eqref{difeq2} $T'_1(x)\to -\infty$ при $x\to +0$, поэтому
из неравенств \eqref{difeq7} и начальных условий \eqref{difeq1} следует,
что $\Ds\varliminf_{x\to +0}T'_k(x)=-\infty$ при всех $k\in \N$.
Таким образом, система уравнений \eqref{difeq} имеет особенность в точке
$x=0$, и поэтому непосредственно решить задачу Коши для данной системы
с начальными условиями в точке $x=0$ невозможно. Для интегрирования системы
\eqref{difeq} мы используем прием, описанный в доказательстве
леммы~\ref{l:Tan2}, а именно: мы найдем решения $y_k(x)$ и $z_k(x)$ уравнения
\eqref{difeq4} на отрезке $[\delta_k,a] $ с начальными условиями
\begin{equation}
y_k(\delta_k)=T_{k-1}(\delta_k),\quad
z_k(\delta_k)=f_k(\delta_k),
\end{equation}
где $\delta_k>0$ --- малое число, $f_k(x)$ --- неявная функция, определяемая
уравнением \eqref{difeq11}.
Кривые $y_k(x)$ и $z_k(x)$, которые мы будем называть соответственно
$k$-м верхним и $k$-м нижним решениями, ограничивают искомую кривую
сверху и снизу:
$$
  z_k(x)<T_k(x)<y_k(x).
$$
В силу замечания~\ref{Rem1} разность $y_k(x)-z_k(x)$ убывает по $x$ и,
следовательно, имеет место оценка:
\begin{equation}
\label{deltak}
\Delta_k= \max_{x\ge\delta_k} |y_k(x)-z_k(x)|=|y_k(\delta_k)-z_k(\delta_k)|=
T_{k-1}(\delta_k)-f_{k}(\delta_k).
\end{equation}
Из непрерывности функций $T_{k-1}(x)$, $f_k(x)$ и соотношения
$T_{k-1}(0)=f_k(0)=1$  вытекает, что
$T_{k-1}(\delta_k)-f_k(\delta_k)\to 0$ при $\delta_k\to 0$.

Для нахождения функции $T_1(x)$ необходимо протабулировать функцию:
$$
x(T_1)=
-\int\limits_{T_1}^1\frac{dt}{t\ln(1- P_1(t))}.
$$
и найти обратную к ней. Использование табличных значений функции $T_1(x)$
в дальнейших вычислениях нежелательно, так как при вычислении $T_1(x)$
в окрестности точки $x=0$ происходит потеря точности. Правая часть системы
\eqref{difeq} явно от $x$ не зависит, а зависит лишь от $T_k(x)$, поэтому
замена переменных $u=T_1(x)$ позволяет решать последующие уравнения системы
($k\ge 2$), не используя $T_1(x)$. Полагая
$$
\tilde T_k(u)= T_{k+1}\left(T_1^{-1}(u)\right),
 \quad k=1,2\dots, m-1,
$$
имеем
$$
\left.\tilde T'_k(u)= T'_{k+1}(x)/T'_1(x)\right|_{x=T_1^{-1}(u)}.
$$
Поскольку $T'_1(x)\to -\infty$ при $x\to +0$, то переход к
переменной $u=T_1(x)$ приводит к уменьшению абсолютных величин производных
искомых функций, что положительно влияет на точность вычислений. Указанная
замена переменных приводит к системе дифференциальных уравнений на отрезке
$[T_1(a),1]$:
\begin{equation*}
\frac{d\tilde T_k}{du}=
\frac{\tilde\phi_k(u,\tilde T_1,\tilde T_2,\dots,\tilde T_k)}{u\ln p(u)}
\end{equation*}
с начальными условиями
$$
 \tilde T_k(1)=1,  \quad k=1,2\dots, m-1.
$$
Пусть первые $k-1$ функций
$$
 \tilde T_i(u),  \quad i=1,2\dots, k-1
$$
найдены. Тогда $\tilde T_k(u)$ есть решение задачи
\begin{align}
\label{difeqp1}
&\frac{d\tilde T_k}{du}=
\tilde \Phi_k(\tilde T_k,u),\\
&\tilde \Phi_k(t,u)=
\frac{\tilde \phi_k(u,\tilde T_1,\tilde T_2,\dots,\tilde T_{k-1},t)}
{u\ln p(u)}
\end{align}
с начальным условиeм
$$
 \tilde T_k(1)=1.
$$
Начальные условия для верхних и нижних кривых после замены переменных
примут вид:
\begin{align*}
&\tilde y_k(u_k)=\tilde T_{k-1}(u_k);\quad
\tilde z_k(u_k)=\tilde f_k(u_k);\\
&\tilde T_0(u_k)=u_k;\quad k=1,2,\dots,n;
\end{align*}
где $u_k=1-\delta_k$, $\delta_k>0$, а $\tilde f_k(u_k)$ есть решение уравнения
$\tilde\Phi_k(t,u_k)=0$ относительно $t$.
Из \eqref{deltak} в силу строгой монотонности функции $T_1(x)$ следует:
$$
 \max_{u\in[T_1(a),u_k]} |\tilde y_k(u)-\tilde z_k(u)|=
\tilde T_{k-1}(u_k)-\tilde f_k(u_k)\to 0  \text { при } u_k\to 1-0\;(k\ge1).
$$

Опишем кратко численный алгоритм решения игры.
Назначение алгоритма: вычисление цены игры $G_{ak}(P,A)$, где
$P(t)=(P_1(t),t)$, при $k=1,2,\dots,m$ и табулирование функций $T_k(x)$
на отрезке $[a_0,a]$, $a_0>0$ с данным шагом $h$.
Работа алгоритма состоит из 2-х этапов.

1-й этап. Производится вычисление значений функции $T_1(x)$
на отрезке $[a_0,a]$ с шагом $h$ путем решения уравнения
$$ x(T_1)=a_0+(i-1)h;\quad i=1,2\dots, M_a,\;
M_a=\left[\frac{a-a_0}{h}\right]+1,
$$
где $x(t)$ --- функция, определенная формулой
\begin{align}
x(t)=-\int\limits_{t}^1\frac{d\tau}{\tau\ln p(\tau)}.
\end{align}

2-й этап. Производится последовательное вычисление значений
функции $T_k(x)$ на отрезке $[a_0,a]$ с шагом $h$ ($k=2,\dots,m$).
На $k$-м шаге алгоритма при табулировании функции $T_k(x)$
производятся следующие вычисления:
\begin{enumerate}
\item
Табулирование $k$-го верхнего решения $\tilde y_k(u)$ уравнения
\eqref{difeqp1} на отрезке $[T_1(a),u_0]$ ($u_0<1$) с начальным условием
$\tilde y_k(u_0)=\tilde T_{k-1}(u_0)$.
\item
Приближенное решение уравнения $\tilde\Phi_k(t,u_0)=0$ относительно $t$
на отрезке $[0,u_0]$. Найденный корень уравнения обозначим $\tilde f_k$.
\item
Табулирование $k$-го нижнего решения $\tilde z_k(u)$ уравнения
\eqref{difeqp1} на отрезке $[T_1(a),u_0]$  с начальным условием
$\tilde z_k(u_0)=\tilde f_k$.
\item
Табулирование функции $\tilde T_k(u)$ по формуле
$$
\tilde T_k(u)=(\tilde y_k(u)+\tilde z_k(u))/2
$$
на отрезке [$T_1(a),u^*]$, где
$$u^*=\Ds\max\{u<u_0:\tilde y_k(u)-\tilde z_k(u)<\varepsilon\},$$
$\varepsilon$ --- заданная точность.
\item
Возвращение к исходной переменной $x$, табулирование функции $T_{k}(x)$
на отрезке $[a_0,a]$.
\item
Вычисление цены игры $v_k(a)$ по формуле \eqref{fin}.
\end{enumerate}